\documentclass{article}

\usepackage{arxiv}

\usepackage[utf8]{inputenc} 
\usepackage[T1]{fontenc}    
\usepackage{hyperref}       
\usepackage{url}            
\usepackage{booktabs}       
\usepackage{nicefrac}       
\usepackage{microtype}      
\usepackage{lipsum}

\usepackage{cite}

\usepackage{amsthm}

\usepackage{amsfonts}
\usepackage{amssymb,amsmath,latexsym,tensor}
\usepackage{mathrsfs}
\usepackage{enumerate}
\usepackage{savesym}
\usepackage{graphicx}
\usepackage[pdftex]{color}
\usepackage[font=footnotesize,labelfont=bf]{caption}
\usepackage[font=footnotesize,labelfont=bf]{subcaption}
\savesymbol{AND}
\usepackage{epstopdf}
\usepackage{algorithm}
\usepackage{algorithmic}
\newcommand{\N}{\mathbb{N}}
\newcommand{\Z}{\mathbb{Z}}

\newcommand{\R}{\mathbb{R}}
\newcommand{\C}{\mathbb{C}}

\newmuskip\pFqmuskip

\newcommand*\pFq[6][8]{%
  \begingroup 
  \pFqmuskip=#1mu\relax
  \mathcode`\,=\string"8000
  \begingroup\lccode`\~=`\,
  \lowercase{\endgroup\let~}\pFqcomma
  {}_{#2}F_{#3}{\left[\genfrac..{0pt}{}{#4}{#5};#6\right]}%
  \endgroup
}
\newcommand{\pFqcomma}{\mskip\pFqmuskip}

\usepackage{mathtools}

\newtheorem{theorem}{Theorem}[section]

\newtheorem{remark}[theorem]{Remark}

\usepackage{mathtools}

\DeclarePairedDelimiter\floor{\lfloor}{\rfloor}

\usepackage{mathtools}
\DeclarePairedDelimiterX{\norm}[1]{\lVert}{\rVert}{#1}
\DeclarePairedDelimiterX{\inp}[2]{\langle}{\rangle}{#1, #2}

\title{Logarithmic Integrals: A review from Gradshteyn and Ryzhik to recent times}

\author{
 Md Sarowar Morshed \\
  Mechanical \& Industrial Engineering\\
  Northeastern University\\
  360 Huntington Ave, Boston, MA 02115\\
  \texttt{morshed.m@husky.neu.edu}
}

\begin{document}
\maketitle
\begin{abstract}

The need to evaluate Logarithmic integrals is ubiquitous in essentially all quantitative areas including mathematical sciences, physical sciences. Some recent developments in Physics namely Feynman diagrams deals with the evaluation of complicated integrals involving logarithmic functions. This work deals with a systematic review of logarithmic integrals starting from \textit{Malmsten} integrals to classical collection of Integrals, Series and Products by \textit{I. S. Gradshteyn} and \textit{I. M. Ryzhik} \cite{ryzhik} to recent times. The evaluation of these types of integrals involves higher transcendental functions (i.e., Hurwitz Zeta function, Polylogarithms, Lerch Transcendental, Orthogonal Polynomials, PolyGamma functions). In a more general sense the following types of integrals are considered for this work:
\begin{align*}
    \int_{0}^{a} f(x) \ln{\{g(x)\}} \ dx
\end{align*}
with $0 \leq a \in \R$ , $f(x)$ and $g(x)$ both either rational/trigonometric or both type of functions.
\end{abstract}


\section{Introduction}
This work devoted to a systematic collection and analysis of logarithmic integrals arises in mathematical literature and from complicated physical models, from the discovery of calculus to the explosion of research in the eighteenth century in evaluating complicated integrals arising from analytic number theory. The most famous example of important logarithmic integrals one can give is the following simple and elegant integral:
\begin{align}
\label{eq:euler}
    \gamma = \lim_{n \xrightarrow{}\infty} \left(-\ln n - H_n^{(1)}\right) = \int_{0}^{\infty} e^{-x} \ln x \ dx
\end{align}
where $\gamma$ is the so called Euler-Mascheroni constant. By using successive differentiation, one can immediately get the following generalized identity:
\begin{align}
\label{eeq:eulergen}
 \int_{0}^{\infty} e^{-\mu x} x^{n} \ln{x} \ dx = \frac{n!}{\mu^{n+1}} \left[H_n-\gamma - \ln{\mu}\right]
\end{align}
where, $n \in \N$ and $H_n$ are the so called Harmonic numbers. Due to the vastness of available logarithmic integrals types, it will be very hard though not impossible to do a systematic review of all kind of logarithmic integrals. Motivated by the richness, complexity and historical importance we devoted our work to the following types of integrals:
\begin{align}
\label{def:li}
    \int_{0}^{a} f(x) \ln{\{g(x)\}} \ dx
\end{align}
with $0 \leq a \in \R$ , $f(x)$ and $g(x)$ both either rational/trigonometric or both type of functions. For example, with the choices $f(x) = e^{-x}, \ g(x)= x$ and $f(x) = e^{-\mu x }x^n, \ g(x)= x$ we can recover the identities shown in equations \eqref{eq:euler} and \eqref{eeq:eulergen} respectively. Moreover, this types of integrals are considered recently by \cite{vardi:1988,Adamchik:1997,boros:2006,moll:2008}.

\subsection{Notation \& Preliminaries}
\label{sec:not}

We follow the standard notation in this work. For example, $\R$, $\N$ and $\Z$ will be used to denote the set of real numbers, the set of integers and the set of natural numbers respectively. $\Re(z)$ is defined as the real part of any complex number $z$. Furthermore,  in evaluating complicated logarithmic integrals we need to use some higher transcendental functions. The Riemann Zeta function and Hurwitz Zeta function are defined as,
\begin{align}
\label{def:zeta}
   \zeta(s)=\sum\limits_{k=1}^{\infty}\frac{1}{k^s}\qquad \text{and} \qquad  \zeta(s,a+1)=\sum\limits_{k=1}^{\infty}\frac{1}{(k+a)^s}=\zeta(s)-H_a^{(s)} 
\end{align}
where $H_n^{(s)}$ the generalized harmonic numbers defined by,
\begin{align}
\label{def:har}
  H_n^{(s)}=\sum\limits_{k=1}^{n}\frac{1}{k^s}=\frac{(-1)^s}{\Gamma (s)}\left[\psi^{(s-1)}(1)-\psi^{(s-1)}(n+1)\right]  
\end{align}
for $ s\geq 1$ and $ n \in \C$. In the above definition, we used the $\psi$ notation for denoting the generalized Poly-gamma function of order $m$ which is given by the $m+1$ th times logarithmic derivative of Gamma function
\begin{align}
\label{def:polgam}
  \psi^{(m)}(z)=\frac{d^m}{dz^m}\psi^{(1)}(z)=\frac{d^{m+1}}{dz^{m+1}}\ln \Gamma (z)  
\end{align}
The Gamma \& Beta functions are defined by 
\begin{align}
\label{def:gamma}
    \Gamma (z) = \int_{0}^{\infty} e^{-x} x^{z-1} dx  \quad \text{and} \quad \mathcal{B} (m,n) = \int_{0}^{1}  x^{m-1} (1-x)^{n-1} dx = \frac{\Gamma (m)\Gamma (n)}{\Gamma (m+n)}
\end{align}

The famous Euler-Mascheroni constant denoted as $\gamma$ can be defined as
\begin{align}
\label{def:eulmas}
   \gamma = \lim_{n \xrightarrow{}\infty} \left(-\ln n - H_n^{(1)}\right) = \int_{0}^{\infty} e^{-x} \ln x \ dx =  -\psi^{(0)}(1) 
\end{align}
The derivatives of generalized harmonic number are given by
\begin{align*}
    \frac{d^m}{dn^m}H_n^{(p)}=(-1)^{m+1}\frac{\Gamma(m+p)}{\Gamma(p)} \zeta(m+p,n+1)
\end{align*}

The Lerch Zeta function, a generalization of zeta functions and its derivative for $\ |z| \leq 1, \ \Re(s) >1$, are defined by
\begin{align}
\label{def:lerch}
    \Phi(z,s,u) = \sum\limits_{k=0}^{\infty} \frac{z^k}{(k+u)^s}, \quad \Phi^{(1)}(z,s,u) = \frac{\partial}{\partial s} \Phi(z,s,u) = -\sum\limits_{k=0}^{\infty} \frac{z^k \ln{(k+u)} }{(k+u)^s}, \quad u > 0, 
\end{align}

The Poly-Logarithm function and its derivative are defined as
\begin{align}
\label{def:polylog}
    \text{Li}_{m} (x) = \text{PolyLog} [m,x] = \sum\limits_{k=1}^{\infty} \frac{x^k}{k^m}, \quad \text{Li}_{m}^{(1)} (x) = \text{PolyLog}^{(1)} [m,x] = \frac{\partial}{\partial m}\text{Li}_{m} (x) =   -\sum\limits_{k=1}^{\infty} \frac{x^k \ln{k}}{k^m}
\end{align}

And we used $(\lambda)_v$  for the Pochhammer symbol defined (for $\lambda,v \in C$ and in terms of the Gamma
function) by
\[(\lambda)_v=\frac{\Gamma(\lambda+v)}{\Gamma(\lambda)}=\lambda(\lambda+1)\cdots (\lambda+v-1); \qquad (\lambda)_0=1\]

The exponential integral $\text{Ei} (x)$ is defined for $x < 0$ and the Catalan constant $G$ defined as follows
\begin{align*}
    \text{Ei} (x) = - \int_{-x}^{\infty} \frac{e^{-t}}{t} dt, \quad \quad  G = \sum\limits_{k=0}^{\infty} \frac{(-1)^k}{(2k+1)^2}
\end{align*}

$B_{k}(x)$, the Bernoulli Polynomials of order $k$ has the following generating function:
\begin{align}
\label{def:bern}
    \frac{t e ^{xt}}{e^t-1} = \sum\limits_{k=0}^{\infty} B_k(x) \frac{t^k}{k!}
\end{align}
The Eulerian polynomials $E_{m}(x)$ are defined by the following generating function:
\begin{align}
\label{def:eulpol}
  \frac{1-x}{1-x e^{t(1-x)}}  = \sum\limits_{k=0}^{\infty} E_k(x) \frac{t^k}{k!}
\end{align}
where, $E_k(x)$ has the following expressions in terms of Eulerian numbers and poly-logarithms:
\begin{align}
\label{eq:eul1}
    E_k(x) = \sum \limits_{j=0}^{k} E_{k,j} x^j, \qquad E_k(-x) = (x+1)^{k+1} \ \text{Li}_{-k}(-x)
\end{align}
where $k \in \N$ and $E_{k,j}$ are the Eulerian numbers. They have the following closed form identity and recurrence formula:
\begin{align}
\label{eq:eul2}
    E_{k,j} = \sum \limits_{i=0}^{j} (-1)^i \binom{k+1}{i} (j-i)^k, \qquad E_{k,j} = j E_{k-1,j} + (k-j+1) E_{k-1,j-1}
\end{align}
Now, let us define polynomials $T_k(x)$ for all $k \in \N$ as follows:
\begin{align}
\label{eq:T1}
    T_k(x) = - \frac{E_{k+1}(-x)}{x} =  \sum \limits_{j=0}^{k} (-1)^j E_{k+1,j+1} x^j
\end{align}
The complete Bell polynomials of order $n$ are defined as
\begin{align}
\label{def:bell}
   Y_n\left[x_1,x_2 \cdots x_n\right]=\sum\limits_{k_1+2k_2+3k_3+\cdots+ nk_n=n}\frac{n!}{k_1!k_2!\cdots k_n!}\left(\frac{x_1}{1!}\right)^{k_1}\left(\frac{x_2}{2!}\right)^{k_2}\cdots \left(\frac{x_n}{n!}\right)^{k_n} 
\end{align}

Let us now consider a function $f(x)$ which has a Taylor series expansion around $x$; then from \cite{ripon:2015}, we have
\begin{align}
    \frac{d^m}{dx^m} e^{f(x)} = e^{f(x)} Y_m\left[f^{1}(x),f^{2}(x) \cdots f^{m}(x)\right] 
\end{align}

The polynomials $P^l(m,x)$ defined as follows:
\begin{align}
\label{def:polynom1}
 P^l(m,x) & = \binom{m}{l} x^l -   \binom{m}{l} \frac{l(m+1)}{2} x^{l-1} + \sum\limits_{k=2}^{l} (-1)^k S(m+1,m+1-k) x^{l-k} \nonumber \\
 & = \binom{m}{l} x^l -   \binom{m}{l} \frac{l(m+1)}{2} x^{l-1} + \sum\limits_{k=2}^{l}  s(m+1,m+1-k) x^{l-k}
\end{align}
where $S(m, k)$ and $s(m, k)$, respectively, denotes the unsigned and signed Stirling numbers of the first kind. Also the sequence $\prescript{r}{k}{\mathbf{C}} (n)$ is defined as to be the numbers are triangle of refined rencontres numbers. The closed form expression for $\prescript{r}{k}{\mathbf{C}} (n)$ can be derived by means of Bell polynomials $Y_k$ defined in equation \eqref{def:bell}:
\begin{align}
\label{polynom2}
    \prescript{r}{k}{\mathbf{C}} (n) = Y_k \left[\gamma + \ln{2n}, \zeta(2), 2! \zeta(3), \cdots, (k-1)\zeta(k)\right]
\end{align}
Note that, the polynomials $P^l(m,x)$ have the following generating function.
\begin{align*}
    (k+1)_m = \sum\limits_{l=0}^{m} (-1)^l P^l\left(m,\frac{p}{n}\right) \left(k+\frac{p}{n}\right)^{m-l}
\end{align*}
A detailed discussion about the analysis \& usefulness of the polynomials $P^l(m,x)$ and $\prescript{r}{k}{\mathbf{C}} (n)$ can be found in \cite{ripon:2015,ripon:2014,ripon:2016}. Finally, we will define some revised notations $\zeta^k(z,q)$ and $\overline{\zeta^k}(z,q) $ for specific types of functions involving zeta functions. The revised zeta functions $\overline{\zeta}(z,q) $ can be defined as follows:
\begin{align}
\label{def:morshed1}
    & \overline{\zeta}(z,q) = \zeta(z,q) - \zeta\left(z, q+\frac{1}{2}\right), \quad \overline{\zeta}(1,x) = \psi \left(x+\frac{1}{2}\right) - \psi (x) 
\end{align}

When, $z = -n$, $\ n \in \N$ we get the following special formula:
\begin{align}
\label{def:morshed2}
 \overline{\zeta}(-n,x) = \frac{1}{n+1} \left[B_{n+1}\left(x+\frac{1}{2}\right)-B_{n+1}(x)\right]  
\end{align}
The functions $\zeta^k(z,q)$ and $\overline{\zeta^k}(z,q) $ are defined as follows:
\begin{align}
    \label{def:morshed4}
    \zeta^k(z,q) =  \frac{\partial^k}{\partial s^k} \zeta(s,q) \bigg|_{s=z},    \quad \overline{\zeta^k}(z,q) =  \frac{\partial^k}{\partial s^k} \overline{\zeta}(s,q) \bigg|_{s=z}
\end{align}
The special case of second function, when $k=1$, we get
\begin{align}
\label{def:morshed3}
    \overline{\zeta^1}(1,x) = \ln{\Gamma(x)} - \ln{\Gamma\left(x+\frac{1}{2}\right)}
\end{align}

\section{Main Results}

In this section, we will discuss logarithmic integrals of the form given in equation \eqref{def:li} for different choices of functions $f(x)$ and $g(x)$. Our goal is to discuss some of the important historical results and their generalizations in terms of recent advances. With various choices of the functions, we divide our results into the following subsections.

\subsection{Integrals with \texorpdfstring{$ g(x) = \ln{\left(\frac{1}{x}\right)}$}{TEXT}}


In the following subsection,  we discuss some of the recent advance Logarithmic Integrals containing the function $\ln{\ln{\left(\frac{1}{x}\right)}}$ in the integrand. Most of these integral roots back to Malmsten, was first discovered by Carl Malmsten in 1842. Modern authors refereed them as the so called Vardi’s integral, which is a particular case of a family of integrals generalized by several researchers in recent times. Adamchick considered some classes of logarithmic integrals in his interesting work \cite{Adamchik:1997}. And a class of similar integrals considered by Medina and Moll in their work \cite{moll:2008} can be found
in the classical tables of integrals originated by Gradshteyn and Ryzhik in \cite{ryzhik:7}. The historical connection between these types of specific logarithmic integrals and the Malmsten integrals has been drawn recently in the work \cite{blagouchine:2014}. More recently, a generalization approach has been proposed by \cite{ripon:2015} to generalize most of the integrals arising from the work of Malmsten, Adamchick, Moll \cite{Adamchik:1997,moll:2008}. In \cite{ripon:2015}, the author considered the following types of integrals 
\begin{align}
\label{eq:def}
    I_Q = \int_{0}^{1} Q(x) \ln{\ln{\left(\frac{1}{x}\right)}} dx
\end{align}

with the choice $f(x) = Q(x)$. Whenever the rational function $f(x)$ has cyclotomic polynomials in the denominator, the above integral of \eqref{eq:def} can be always evaluated in terms of higher order derivatives of Lerch zeta function (i.e., see definition \eqref{def:lerch}). However, it is not yet known that for what type of rational functions $f(x)$ the above integral has a closed form solution.
\begin{remark}
Note that, using integration by parts one can derive that the evaluation of integral $I_Q$ can be reduced to evaluation of the following integral
\begin{align}
\label{def:pot}
    J_{Q} = \int_{0}^{1} \frac{Q(x)}{\ln{x}} dx
\end{align}
This class has some important applications in Mathematical physics in the so called Potts model for the triangular lattice as shown by Baxter \textit{et. al} \cite{potts:1978}. From a evaluation standpoint, they are first analyzed by Adamchick in his work \cite{Adamchik:1997}. Adamchick discovered analytic expressions for the integrals given in \eqref{def:pot} with the choice $Q(x)$ as cyclotomic polynomials.
\end{remark}

Now, let us discuss some of the integrals arising in Historical context of the above form (i.e., equation \eqref{eq:def}):

\subsection{\texorpdfstring{$f(x) = x^{m}/ h(x)$, \ $h(x)$}{TEXT} cyclotomic polynomials}

\paragraph{Malmsten's Integrals}

The oldest known integral containing the above form are the integrals considered by Malmsten in 1842. In their work, \cite{malmsten:1,malmsten:2}  Malmsten and his colleagues evaluated the following integral pair:
\begin{align}
   &  \int_{0}^{1} \frac{\ln{\ln{\left(\frac{1}{x}\right)}}}{1+x+x^2}  \ dx = \int_{1}^{\infty} \frac{\ln{\ln{\left(x\right)}}}{1+x+x^2}  \ dx  = \frac{\pi}{\sqrt{3}} \ln \left( \frac{\sqrt[\leftroot{-2}\uproot{2}3]{2\pi} \Gamma(\frac{2}{3})}{\Gamma(\frac{1}{3})}\right) \\
   & \int_{0}^{1} \frac{\ln{\ln{\left(\frac{1}{x}\right)}}}{1-x+x^2}  \ dx = \int_{1}^{\infty} \frac{\ln{\ln{\left(x\right)}}}{1-x+x^2}  \ dx  = \frac{2\pi}{\sqrt{3}} \ln \left( \frac{\sqrt[\leftroot{-2}\uproot{2}6]{32\pi^5}}{\Gamma(\frac{1}{6})}\right)
\end{align}
One can note that, the above rational function in the integrand has  cyclotomic polynomials in the denominator. Malmsten and his colleagues evaluated some other amazing integrals too in their work, but unfortunately later authors missed the reference of Malmsten and renamed these integrals of \eqref{eq:var} as Vardi's integral. According to the recent work \cite{blagouchine:2014}, Malmsten obtained quite a large number of integrals containing the terms $\ln{\ln{\left(\frac{1}{x}\right)}}$ and $\ln{\ln{\left(x\right)}}$ in the numerator. Following are some examples of such integral pair:
\begin{align}
     Y_n & = \int_{0}^{1} \frac{x^{n-2} \ \ln{\ln{\left(\frac{1}{x}\right)}} \ dx}{1+x^2+x^4+...+x^{2n-2}}   = \int_{1}^{\infty} \frac{x^{n-2} \ \ln{\ln{\left(x\right)}} \ dx}{1+x^2+x^4+...+x^{2n-2}}    \label{eq:malm1} \\
   &  = \begin{cases}
        \frac{\pi}{2n} \tan\left(\frac{\pi}{2 \pi}\right) \ \ln{(2\pi)} - \frac{\pi}{n}\sum \limits_{l=1}^{n-1} (-1)^l \ \sin{\left(\frac{\pi l}{n}\right)} \ln \left\{\frac{\Gamma\left(\frac{1}{2}+\frac{l}{2n}\right)}{\Gamma(\frac{l}{2n})}\right\}, \ n = 2,4,6,...\\
        \frac{\pi}{2n} \tan\left(\frac{\pi}{2 \pi}\right) \ \ln{(\pi)} - \frac{\pi}{n}\sum \limits_{l=1}^{\frac{1}{2}(n-1)} (-1)^l \ \sin{\left(\frac{\pi l}{n}\right)} \ln \left\{\frac{\Gamma\left(1-\frac{l}{n}\right)}{\Gamma(\frac{l}{n})}\right\}, \ n = 1,3,5,...
       \end{cases} \nonumber \\
   X_n & = \int_{0}^{1} \frac{x^{n-2} \ \ln{\ln{\left(\frac{1}{x}\right)}} \ dx}{1-x^2+x^4-...+x^{2n-2}}   = \int_{1}^{\infty} \frac{x^{n-2} \ \ln{\ln{\left(x\right)}} \ dx}{1-x^2+x^4-...+x^{2n-2}} \label{eq:malm2} \\
   & = \frac{\pi}{2n} \sec\left(\frac{\pi}{2 \pi}\right) \ \ln{(\pi)} - \frac{\pi}{n}\sum \limits_{l=1}^{\frac{1}{2}(n-1)} (-1)^l \ \cos{\left(\frac{(2l-1)\pi}{2n}\right)} \ln \left\{\frac{\Gamma\left(1-\frac{2l-1}{2n}\right)}{\Gamma(\frac{2l-1}{2n})}\right\} \nonumber
\end{align}

\paragraph{Gradshteyn and Ryzhik Tables} \footnote{this section is under construction} \cite{ryzhik}

\paragraph{Vardi's Integral} Vardi's remarkable paper \cite{vardi:1988} contains an interesting method for calculating integrals defined of the form shown in equation \eqref{eq:def}. In particular Vardi considered the following integral:
\begin{align}
\label{eq:var}
    \int_{0}^{1} \frac{1}{1+x^2} \ln{\ln{\left(\frac{1}{x}\right)}} \ dx = \int_{\pi/4}^{\pi/2}  \ln{\ln{\left(\tan x\right)}} \ dx = \pi \ln \left( \frac{\sqrt{2 \pi}\Gamma(\frac{3}{4})}{\Gamma(\frac{1}{4})}\right)
\end{align}
This integral was first evaluated by Malmsten is his remarkable work \cite{malmsten:1}, which remained unnoticed to the research community until 1998 \cite{vardi:1988}.
The formal method for deriving the  above identity is provided in great detail in \cite{vardi:1988, moll:2008}. Let us define $\Gamma_Q(s)$ as
\begin{align*}
  \Gamma_Q(s) =   \int_{0}^{1} Q(x) \left(\ln{\left(\frac{1}{x}\right)}\right)^{s-1} dx
\end{align*}
One can check that $I_Q = \Gamma_Q^{'}(1)$. Explicit evaluation of $I_Q$ is possible when the function $Q(x)$ is analytic at $x = 0$.

\begin{theorem}
\label{th:vardi}
Assume that the function $Q(x)$ is analytic at $x = 0$. And define the Laurent series expansion of $Q(x)$ and the corresponding $L$- series expansion as follows:
\begin{align*}
    Q(x) = \sum\limits_{k=0}^{\infty} a_k x^k, \quad \text{and} \quad L_Q(s) = \sum\limits_{k=0}^{\infty} \frac{a_k}{(k+1)^s}
\end{align*}
Then, the integrals $I_Q$ can be calculated using the following identity
\begin{align*}
    I_Q = \int_{0}^{1} Q(x) \ln{\ln{\left(\frac{1}{x}\right)}} dx = - \gamma L_Q(1) + L^{'}_{Q}(1)
\end{align*}
\end{theorem}

The above Theorem can be found in \cite{vardi:1988, moll:2008}. Later, we will show how this Theorem can be used to derive some complicated integrals.

\paragraph{Adamchick's Integrals} The following integrals have been considered by Admachick in his work \cite{Adamchik:1997}:
\begin{align}
\int_{0}^{1} \frac{x^{p-1}}{\left(1+x^n\right)^{m+1}} \ln{\ln{\left(\frac{1}{x}\right)}} \ dx, \ m = 0,1,2; \quad    \int_{0}^{1} \frac{x^{p-1}(1-x)}{\left(1-x^n\right)} \ln{\ln{\left(\frac{1}{x}\right)}} \ dx
\end{align}
Notably, this class of integrals found application in statistical physics \& lattice theory. Specifically the Potts model on the triangular lattice \cite{potts:1978,potts:1997} uses these type of integrals. For an elaborate discussion about the technical details we refer to the papers \cite{potts:1978,Adamchik:1997}. Some of the notable integrals from Adamchick's work are stated below:

\begin{theorem}
\label{thm:ad1}
(Proposition 3 in \cite{Adamchik:1997})  Let $\Re(p), \ \Re(n) > 0$, then the following identity holds:
\begin{align*}
    & \int_{0}^{1} \frac{x^{p-1}}{\left(1+x^n\right)} \ln{\ln{\left(\frac{1}{x}\right)}} \ dx = -\frac{\gamma+\ln{2n}}{2n} \ \overline{\zeta}\left(1,\frac{p}{2n}\right) + \frac{1}{2n} \overline{\zeta^{1}} \left(1,\frac{p}{2n}\right)
\end{align*}
\end{theorem}

\begin{theorem}
\label{thm:ad2}
(Proposition 4 in \cite{Adamchik:1997})  Let $\Re(p) >0, \ \Re(n) \geq 1$, then the following identity holds:
\begin{align*}
\int_{0}^{1} \frac{x^{p-1}(1-x)}{\left(1-x^n\right)} & \ln{\ln{\left(\frac{1}{x}\right)}} \ dx \\
& = \frac{\gamma+\ln{n}}{n} \left[\psi\left(\frac{p}{n}\right)-\psi\left(\frac{p+1}{n}\right)\right] + \frac{1}{n} \left[\zeta^{'}\left(1,\frac{p}{n}\right)-\zeta^{'}\left(1,\frac{p+1}{n}\right)\right]
\end{align*}
\end{theorem}

\begin{theorem}
\label{thm:ad3}
(Proposition 5 in \cite{Adamchik:1997})  Let $\Re(p), \ \Re(n) > 0$, then the following identity holds:
\begin{align*}
& \int_{0}^{1} \frac{x^{p-1}}{\left(1+x^n\right)^2} \ln{\ln{\left(\frac{1}{x}\right)}} \ dx = -\frac{(\gamma+\ln{2n})(n-p)}{2n^2} \ \overline{\zeta}\left(1,\frac{p}{2n}\right) + \frac{n-p}{2n^2} \ \overline{\zeta^{1}} \left(1,\frac{p}{2n}\right) \nonumber \\
    &    \qquad \qquad \qquad \qquad \qquad \qquad \qquad -\frac{1}{2n} \left[\gamma+\ln{2n}-2 \ln{\left(\frac{\Gamma(p/2n)}{\Gamma(n+p/2n)}\right)}\right]
\end{align*}
\end{theorem}

\begin{theorem}
\label{thm:ad4}
(Proposition 6 in \cite{Adamchik:1997})  Let $\Re(p), \ \Re(n) > 0$, then the following identity holds:
\begin{align*}
     \int_{0}^{1} \frac{x^{p-1}}{\left(1+x^n\right)^3} \ln{\ln{\left(\frac{1}{x}\right)}} \ dx & = -\frac{(\gamma+\ln{2n})(n-p)(2n-p)}{2n^2} \ \overline{\zeta}\left(-1,\frac{p}{2n}\right) \\
    & + \frac{3n-2p}{2 n^2} \ln{\left(\frac{\Gamma\left(\frac{p}{2n}\right)}{\Gamma\left(\frac{p+n}{2n}\right)}\right)} + \frac{1}{n} \overline{\zeta^{1}} \left(1,\frac{p}{2n}\right) \nonumber \\
    &   + \frac{(2p-5n)(\gamma+\ln{2n})}{8 n^2}+ \frac{(n-p)(2n-p)}{4 n^3} \ \overline{\zeta^{1}} \left(1,\frac{p}{2n}\right)
\end{align*}
\end{theorem}

where the definitions of revised zeta functions $\overline{\zeta}(z,q)$ and $\overline{\zeta^k}(z,q)$ are defined in equation \eqref{def:morshed1} - \eqref{def:morshed4}. For a more detailed discussion about these types of functions see \cite{ripon:2015}.

\paragraph{Morshed's Generalization}

Morshed \cite{ripon:2015} specifically considered the integrals of following types:
\begin{align}
    \int_{0}^{1} \frac{x^{p-1}}{\left(1+x^n\right)^{m+1}} \left[\ln{\ln{\left(\frac{1}{x}\right)}}\right]^r  \ dx \quad \quad \text{and} \quad \int_{0}^{1} \frac{x^{p-1}(1-x)}{\left(1-x^n\right)^{m+1}} \left[\ln{\ln{\left(\frac{1}{x}\right)}}\right]^r  \ dx
\end{align}
where $\Re(p), \Re(n) > 0$ and $m,r \in \N$. By choosing specific set of values for the parameters $p,n,m,r$ one can recover almost all of the integrals appeared in Malmsten, Adamchick and Moll's work. \cite{Adamchik:1997,moll:2008}. The following propositions propoosed by the author in his work \cite{ripon:2015},

\begin{theorem}
Let $\Re(p), \Re(n) > 0$, then for all $m,r \in \N$ the following identity holds:
\begin{align}
    \int_{0}^{1} \frac{x^{p-1}}{\left(1+x^n\right)^{m+1}} \left[\ln{\ln{\left(\frac{1}{x}\right)}}\right]^r  \ dx = \frac{2^{m-1}}{n m!} \sum \limits_{l=0}^{m} \sum \limits_{k=0}^{r} \frac{(-1)^{l+k}}{2^l} \binom{r}{k} \prescript{r}{k}{\mathbf{C}} (n) P^l \left(m,\frac{p}{n}\right) \overline{\zeta^{r-k}}\left(1+l-m, \frac{p}{2n}\right)
\end{align}
\end{theorem}
For example, if we let $r=1, m=0$ we get Proposition 3, $r=1, m=1$ we get Proposition 5 and $r=1, m=2$ we get Proposition 6  given in \cite{Adamchik:1997}.

\begin{theorem}
Let $\Re(p), \Re(n) > 0$, then for all $m,r \in \N$ the following identity holds:
\begin{align}
    \int_{0}^{1} \frac{x^{p-1}(1-x)}{\left(1-x^n\right)^{m+1}}& \left[\ln{\ln{\left(\frac{1}{x}\right)}}\right]^r  \ dx  = \frac{1}{n m!} \sum \limits_{l=0}^{m} \sum \limits_{k=0}^{r} (-1)^{l+k} \binom{r}{k} \prescript{r}{k}{\mathbf{C}} (n) P^l \left(m,\frac{p}{n}\right) \zeta^{r-k}\left(1+l-m, \frac{p}{n}\right) \nonumber \\
    & - \frac{1}{n  m!} \sum \limits_{l=0}^{m} \sum \limits_{k=0}^{r} (-1)^{l+k} \binom{r}{k} \prescript{r}{k}{\mathbf{C}} (n) P^l \left(m,\frac{p+1}{n}\right) \zeta^{r-k}\left(1+l-m, \frac{p+1}{n}\right)
\end{align}
\end{theorem}
For example, if we let $r=1, m=0$ we get Proposition 4 given in \cite{Adamchik:1997}. The polynomials $P^l(m,x)$ and $\prescript{r}{k}{\mathbf{C}} (n)$ are defined in equation \eqref{def:polynom1} and \eqref{polynom2} respectively.

\paragraph{Irresistible Integrals}

The work \cite{boros:2006}, contains the following important result:
\begin{align*}
   \int_{0}^{1} \frac{1}{\left(1-x+x^2\right)^2} \ln{\ln{\left(\frac{1}{x}\right)}} \ dx = -\frac{\gamma}{3} - \frac{1}{3} \ln{\left(\frac{6 \sqrt{3}}{\pi}\right)}+ \frac{\pi \sqrt{3}}{27} \left(5 \ln{(2 \pi)}-6 \ln{\Gamma\left(\frac{1}{6}\right)}\right)
\end{align*}

The authors (\cite{boros:2006} section 9-10) provided some classical identities provided below:
\begin{align*}
    \int_{0}^{\infty} e^{-\mu x} x^{s-1} \ln{x} \ dx = \frac{\Gamma(s)}{\mu^s} \left[\psi(s)- \ln{\mu}\right], \quad  \int_{0}^{\infty} e^{-\mu x} x^{n} \ln{x} \ dx = \frac{n!}{\mu^{n+1}} \left[H_n-\gamma - \ln{\mu}\right]
\end{align*}
where, $s \in \C$ and $n \in \N$. And for $n \in \N$ we have the following identity
\begin{align*}
 \int_{0}^{\infty} e^{-\mu x} x^{n-\frac{1}{2}} \ln{x} \ dx = \frac{\sqrt{\pi} (2n)!}{\mu^{n+\frac{1}{2}} 2^{2n}n!} \left[2 H_{2n}-H_n-\gamma - \ln{4\mu}\right]
\end{align*}

And the following identity involving beta function
\begin{align*}
    \int_{0}^{1} x^{as-1} \left(1-x^a\right)^{z-1} \ln x \ dx = \frac{\psi(s)\Gamma(z)\Gamma(s+z)-\psi(s+z)\Gamma(s)}{a^2 \Gamma^2(s+z)}
\end{align*}

The following identities can be derived by some elementary operations:
\begin{align*}
    & \int_{0}^{1} \frac{\ln^k{x}}{1+x} \ dx = \frac{(-1)^k k! \left(2^k-1\right)}{2^k} \ \zeta(k+1), \qquad \int_{0}^{1} \frac{\ln^k{x}}{1-x^2} \ dx = \frac{(-1)^k k! \left(2^{k+1}-1\right)}{2^{k+1}} \ \zeta(k+1) \\
    & \int_{0}^{1} \frac{\ln{(1-x^2)}}{x} \ln^{k-1}{x} \ dx = \frac{(-1)^k (k-1)!}{2^{k}} \ \zeta(k+1), \qquad \int_{1}^{p} \frac{\ln{(1-x)}}{x} \ dx = \frac{\pi^2}{6} - \text{Li}_{2}(p)
\end{align*}

The following logarithmic integrals can be derived as closed form expressions with Bernoulli \& Euler polynomials:
\begin{align*}
    \int_{0}^{1} \frac{\ln^{2n-1}x}{1-x^2} \ dx = \frac{\pi^{2n}\left(2^{2n}-1\right)}{4n} B_{2n}, \quad \int_{0}^{1} \frac{\ln^{2n+1}x}{1+x^2} \ dx = \frac{\pi^{2n+1}}{2^{2n+2}} E_{2n}
\end{align*}
for $n \in \N$. 
The above transformation can be proven with some basic elementary manipulations using substitution. With more careful analysis one can derive the following complicated integrals for $k \in \N$:
\begin{align*}
    & \int_{0}^{1} \ln{(1+x)}\ln^{k-1} x \ \frac{ dx }{x}  =  (-1)^{k-1} \ (k-1)! \left(1-2^{-k}\right) \ \zeta(k+1) \\
    & \int_{0}^{1} \ln{(1-x)}\ln^{k-1} x \ \frac{ dx }{x}  =  (-1)^{k-1} \ (k-1)!\ \zeta(k+1)
\end{align*}

Note that, most of the above identities can also be found on the classical tables of Gradshteyn \& Ryzhik \cite{ryzhik}. In \cite{boros:2006}, (section 5) authors provided the following beautiful integral sequence. Let's define $f_n (x)$ for $n \geq 0, n \in \N, a >0$ such that
\begin{align}
    f_0(x) = \ln(a+x), \quad f_n(x) = \int_{0}^{x} f_{n-1}(t) \ dt
\end{align}
Then, $f_n(x)$ has the following closed form solution:
\begin{align*}
    f_n(x) = \frac{(x+a)^n}{n!} \left[\ln(x+a) - H_n\right] + \ln a \left[x^n - (x+a)^n\right] + \frac{1}{n!} \sum\limits_{k=1}^{n} \binom{n}{r} a^k x^{n-k} H_k
\end{align*}

The following interesting generating function  can be used to derive many complicated integrals 
\begin{align*}
    \ln{(1+x)} \ln{(1-x)} = \sum\limits_{k=1}^{\infty} \left(H_k-H_{2k}-\frac{1}{2}\right) \frac{x^{2k}}{k}
\end{align*}
like the following:
\begin{align*}
    & \int_{0}^{1} \ln{(1+x)} \ln{(1-x)} \ dx = \ln^2{2} - 2 \ln{2} +2 - \zeta(2)
\end{align*}
The following closed form evaluation obtained explicitly:
\begin{align*}
    \int_{a}^{b} \frac{\ln{x} \ dx}{q_1 x+q_0} = \frac{1}{q_1} \left\{\ln{b} \ \ln{(1+bq^*)}-\ln{a} \ \ln{(1+aq^*)}- \text{Li}_2[-aq^*]+ \text{Li}_2[-bq^*]\right\}
\end{align*}
with the choice $q^* = q_1/q_0$. The following expressions widely known as Neilsen-Ramanujan constants
\begin{align*}
    a_k = \int_{1}^{2} \frac{\ln^k x}{x-1} \ dx, \quad \text{for} \ k \geq 1
\end{align*}
These were studied by Gosper (1996) and the following closed form relation was derived for the sequence $a_k$,
\begin{align}
\label{eq:gosper}
    a_k = k! \ \zeta (k+1) - \frac{k}{k+1} \ln^{k+1} 2 - k! \sum\limits_{j=0}^{k-1} \frac{\ln^j 2}{j!} \text{Li}_{k+1-j} \left(\frac{1}{2}\right)
\end{align}
with the special values $a_1 = \zeta(2)/2, \ a_2 = \zeta(3)/4$. It can be shown that the following  generalized Theorem holds
\begin{theorem}
\label{th:gen}
For all $1 \leq b \in R^+$, $ 1 \leq k \in \Z$ and $r \neq 1,2,3,...$, we have
\begin{align*}
    & a_k(b,r) = \int_{1}^{b} \frac{x^r \ \ln^k x}{x-1} \ dx =  k! \ \zeta (k+1,1-r) + g(r,k,b)   - \sum\limits_{j=0}^{k} \frac{ k!}{j!} \left(\ln{b}\right)^{j} b^{r-1} \Phi\left(\frac{1}{b},k-j+1,1-r\right) \\
    & g(0,k,b) = \frac{\left(\ln{b}\right)^{k+1}}{k+1}, \qquad g(r,k,b) =  \int_{0}^{\ln{b}} u^k e^{ur} \ du = \sum\limits_{j=0}^{k} \frac{(-1)^{k+j}}{r^{k-j+1}}\frac{ k!}{j!} \left(\ln{b}\right)^{j}
\end{align*}
\end{theorem}
with the choice of $b =2, r = 0$, after some simplification we can get the special identity in equation \eqref{eq:gosper}.
\begin{proof}
Note that, with some elementary substitution we get the following equivalent identity for $a_k(b)$:
\begin{align}
\label{eq:pr1}
    a_k(b) = \int_{1}^{b} \frac{ x^r \ln^k x}{x-1} \ dx & =  \int_{0}^{\ln{b}} u^k e^{ur} \ du + \int_{0}^{\infty} \frac{u^k e^{ur} \ du}{e^u -1} - \int_{\ln{b}}^{\infty} \frac{u^k e^{ur} \ du}{e^u -1} \nonumber \\
    & = g(r,k,b) + k! \ \zeta(k+1,1-r)- \int_{\ln{b}}^{\infty} \frac{u^k e^{ur} \ du}{e^u -1} 
\end{align}
where we used the following well-known identity for zeta function:
\begin{align*}
    \Gamma(k+1) \zeta(k+1,1-r) = \int_{0}^{\infty} \frac{u^k e^{ur} \ du}{e^u -1}
\end{align*}
Now we can simply evaluate the third part of the equation \eqref{eq:pr1} as follows
\begin{align}
\label{eq:pr2}
    \int_{\ln{b}}^{\infty} \frac{u^k e^{ur} \ du}{e^u -1}  = \sum\limits_{n=1}^{\infty}  \int_{\ln{b}}^{\infty} u^k e^{ur-un} \ du = \sum\limits_{n=1}^{\infty} \frac{1}{(n-r)^{k+1}}  \int_{(n-r)\ln{b}}^{\infty} t^k e^{-t} \ dt 
\end{align}
Now, note that, the definite integral can be calculated as follows:
\begin{align}
\label{eq:pr3}
   \int_{(n-r)\ln{b}}^{\infty} t^k e^{-t} \ dt  = - \sum\limits_{j=0}^{k} \frac{k!}{j!} \left( e^{-t} \ t^j\right) \bigg|_{x= (n-r) \ln{b}}^{x= \infty} = \sum\limits_{j=0}^{k} \frac{k!}{j!} \frac{(n-r)^j}{b^{n-r}} \left(\ln{b}\right)^j
\end{align}
Now, using the integral identity of equation \eqref{eq:pr3} to equation \eqref{eq:pr2}, we get the following:
\begin{align}
\label{eq:pr4}
     \int_{\ln{b}}^{\infty} \frac{u^k e^{ur} \ du}{e^u -1} & = \sum\limits_{n=1}^{\infty} \frac{1}{(n-r)^{k+1}}  \int_{(n-r)\ln{b}}^{\infty} t^k e^{-t} \ dt  \\
    & = \sum\limits_{j=0}^{k} \frac{k!}{j!} \left(\ln{b}\right)^j  \sum\limits_{n=0}^{\infty} \frac{b^{r-1}\left(\frac{1}{b}\right)^n}{(n-r+1)^{k-j+1}}  = \sum\limits_{j=0}^{k} \frac{ k!}{j!} \left(\ln{b}\right)^{j} \Phi\left(\frac{1}{b},k-j+1,1-r\right)
\end{align}
Now, we have
\begin{align*}
    g(0,k,b) = \int_{0}^{\ln{b}} u^k \ du = \frac{\left(\ln{b}\right)^{k+1}}{k+1}
\end{align*}
And for $r \neq 0,1,2,3,,,,$ we have,
\begin{align*}
    g(r,k,b) = \int_{0}^{\ln{b}} u^k e^{ur} \ du = \left(-\frac{1}{r}\right)^{k+1}  \int_{0}^{-r\ln{b}} y^k e^{-y} \ dy = \sum\limits_{j=0}^{k} \frac{(-1)^{k+j}}{r^{k-j+1}}\frac{ k!}{j!} \left(\ln{b}\right)^{j}
\end{align*}
where we used the idea from equation \eqref{eq:pr3}. Now combining equation \eqref{eq:pr1} , \eqref{eq:pr4} and the derivation of $g(r,k,b)$, we get the desired result of Theorem \ref{th:gen}.
\end{proof}

They also proposed the following vanishing integral formula for any $a \in \R^+$:
\begin{align*}
    \int_{0}^{\infty} \frac{x^{a-1} \ln^n x}{\left(1+x^{2a}\right)^{n+1}} C_n \left(-x^{2a}\right) \ dx = 0
\end{align*}
where the polynomials $C_n(x)$ satisfies the following recurrence relation:
\begin{align*}
    C_{n+1}(x) = 2 x(1-x) C_n^{'}(x) + [1+x+2nx] C_n(x)
\end{align*}
It can be noted that $C_n(x)$ has the following closed form:
\begin{align*}
    C_n(x) = \sum\limits_{k= 0}^{n} x^{n-k} T(n,k)
\end{align*}
where, $T(n,k)$ are the sol called triangle numbers.

The evaluation of following integral family widely known as logsin function can be found in their work:
\begin{align*}
  S_n & = \int_{0}^{\pi} \left(\ln{(\sin{x})}\right)^n \ dx  = 2^n \ \int_{0}^{\frac{\pi}{2}} \left(\ln{(\sin{x})}\right)^n \ dx \\
  & = 2^n \ \int_{0}^{\frac{\pi}{2}} \left(\ln{(\cos{x})}\right)^n \ dx =  \frac{(-1)^{n-1}}{(n-1)!} \int_{0}^{\frac{\pi}{2}} \left(\sin{\theta}\right)^{n-1} 
\end{align*}
They proved the following recurrence relations for the function $S_n$
\begin{align*}
   &  S_1 S_{2n-1} - S_2 S_{2n-2} +...+ S_{2n-1} S_1 = (-1)^{n-1} \frac{\left(2^{2n}-1\right)}{(2n)!} \pi^{2n} B_n \\
   & (n-1)S_{n} = \ln{2} \ S_{n-1} + \sum\limits_{k=1}^{n-1} \left(1-2^{-k}\right) \zeta_{k+1} S_{n-k-1}  
\end{align*}

The following closed form expression can be found for an $n > 1$:
\begin{align*}
    \int_{0}^{1} \frac{\ln{x}}{ \sqrt[\leftroot{-2}\uproot{2}n]{1-x^{2n}}} \ dx = - \frac{\pi \mathcal{B}\left(\frac{1}{2n},\frac{1}{2n}\right)}{8 n^2 \sin{\left(\frac{\pi}{2n}\right)}}, \qquad \int_{0}^{1} \frac{\ln{x}}{ \sqrt[\leftroot{-2}\uproot{2}n]{x^{n-1}(1-x^{2n})}} \ dx = - \frac{\pi \mathcal{B}\left(\frac{1}{2n},\frac{1}{2n}\right)}{8  \sin{\left(\frac{\pi}{2n}\right)}}
\end{align*}

The following master integral has been used in \cite{boros:2006}, to derive some complicated integrals:
\begin{align}
\label{eq:bor10}
    g(r,a) = \int_{0}^{\infty} \left[\frac{x^2}{x^4+2ax^2+1}\right]^r \frac{x^2+1}{x^2(x^s+1)} \ dx = \frac{\mathcal{B}\left(r-\frac{1}{2},\frac{1}{2}\right)}{2^{r+\frac{1}{2}}(1+a)^{r+\frac{1}{2}}}
\end{align}
One can note that for some fixed values of $a$, differentiating the function $g(r,a)$, with respect to $r$, the following integrals can be derived easily:
\begin{align*}
    & \int_{0}^{\infty} \left(\frac{x^2}{x^4-x^2+1}\right)^{\frac{3}{4}} \ln{\left(\frac{x^2}{x^4-x^2+1}\right)} \ dx = - \frac{\sqrt{\pi}}{2\sqrt{2}} \ \Gamma^2\left(\frac{3}{4}\right) \\
    & \int_{0}^{\infty} \left(\frac{1}{x^4-x^2+1}\right) \ln^2{\left(\frac{x^2}{x^4-x^2+1}\right)} \ dx = \frac{\pi}{2} \left(\frac{\pi^2}{3} +4 \ln^2{2}\right)
\end{align*}

With defining the following function
\begin{align*}
    G(r) = \frac{2}{\sqrt{\pi}} \int_{0}^{\infty} \left(\frac{x}{x^2+1}\right)^{2r} \ \frac{dx}{x^2} = \frac{\Gamma\left(r-\frac{1}{2}\right)}{\Gamma (r) \ 2^{2r-1}}
\end{align*}
Differentiating $G(r)$ $n$ times with respect to $r$,  we get the following master formula:
\begin{align*}
   \int_{0}^{\infty} \left(\frac{x}{x^2+1}\right)^{2r}  \left[\ln {\left(\frac{x}{x^2+1}\right)}\right]^n \ \frac{dx}{x^2} & = \frac{\sqrt{\pi}}{2^{n}} \ \frac{d^{n}}{dr^{n}} \left[\frac{\Gamma\left(r-\frac{1}{2}\right)}{\Gamma (r) \ 2^{2r-1}}\right] \\
   & = \frac{\left(\ln{2}\right)^n}{2^{2r}} \sum\limits_{k=0}^{n} \frac{(-1)^{n-k}}{\left(2\ln{2}\right)^k} \binom{n}{k} \ \frac{d^{k}}{dr^{k}} \left\{ \mathcal{B} \left(r-\frac{1}{2},\frac{1}{2}\right)\right\}
\end{align*}
With the choice $r = \frac{3}{2}$ and $n=4$, we get the following result:
\begin{align*}
   \int_{0}^{\infty} \left(\frac{x}{x^2+1}\right)^{3}  \left[\ln {\left(\frac{x}{x^2+1}\right)}\right]^4 \ \frac{dx}{x^2}  = - \frac{3 \pi^4}{160}- \frac{\pi^2}{2} + 12 - 3 \zeta(3)
\end{align*}
With $r = 2m+1$ and $n = 1$, we get the following closed form in terms of Harmonic numbers:
\begin{align*}
   \int_{0}^{\infty} \left(\frac{x}{x^2+1}\right)^{2m+1}  \ln {\left(\frac{x}{x^2+1}\right)}\ \frac{dx}{x^2}  = \frac{H_m-H_{2m}-\frac{1}{2m}}{2m \binom{2m}{m}}
\end{align*}

With the help of some elementary operations one can evaluate the following integral:
\begin{align*}
    \int_{0}^{\infty} \left(\frac{x}{x^2+1}\right)^{2r+1} \left\{\frac{\ln{\left(1-x+x^2-,,,+x^{2p}\right)}}{\ln{x}}\right\} \ \frac{dx}{x} = \frac{\pi r p \ \Gamma(2r)}{2^{4r} \ \Gamma^2(r+1)}
\end{align*}
In their book one can also find the following vanishing integral:
\begin{align*}
 \int_{0}^{\infty} \frac{x^{b+c-1} \ \ln{x}}{\left(1+x^2\right)^c \left(1+x^b\right)^2} \ dx = 0, \qquad \int_{0}^{\infty} \frac{x^{b+c-1} (x^b-1) \ \ln^2{x}}{\left(1+x^2\right)^c \left(1+x^b\right)^3} \ dx = 0
\end{align*}

\paragraph{Malmsten Integrals} In \cite{malmsten:1}, using some lengthy calculation, the following integrals are evaluated:
\begin{align*}
     \int_{0}^{1} \frac{\sinh{(bx)}}{\sinh{(\pi x)}} & \ln{(a^2+x^2)} \ dx \\
     & = \begin{cases}
         \tan\left(\frac{\pi m}{2 n}\right) \ \ln{(2n)} -2 \sum \limits_{l=1}^{n-1} (-1)^l \ \sin{\left(\frac{\pi ml}{n}\right)} \ln \left\{\frac{\Gamma\left(\frac{1}{2}+\frac{l+a}{2n}\right)}{\Gamma(\frac{l+a}{2n})}\right\}, \ m+n \ \text{is odd} \\
        \tan\left(\frac{\pi m}{2 n}\right) \ \ln{(n)} -2 \sum \limits_{l=1}^{\floor {\frac{1}{2}(n-1)}} (-1)^l \ \sin{\left(\frac{\pi ml}{n}\right)} \ln \left\{\frac{\Gamma\left(1-\frac{l-a}{n}\right)}{\Gamma(\frac{l+a}{n})}\right\}, \ m+n \ \text{is even}
       \end{cases} \\
   \int_{0}^{1} \frac{\cosh{(bx)}}{\cosh{(\pi x)}} &  \ln{(a^2+x^2)} \ dx  \\
   & = \begin{cases}
         \sec\left(\frac{\pi m}{2 n}\right) \ \ln{(2n)} -2 \sum \limits_{l=1}^{n-1} (-1)^l \ \cos{\left(\frac{(2l-1)m \pi}{n}\right)} \ln \left\{\frac{\Gamma\left(\frac{1}{2}+\frac{2l+2a-1}{4n}\right)}{\Gamma(\frac{2l+2a-1}{4n})}\right\}, \ m+n \ \text{is odd} \\
        \sec\left(\frac{\pi m}{2 n}\right) \ \ln{(n)} -2 \sum \limits_{l=1}^{\floor {\frac{1}{2}(n-1)}} (-1)^l \ \cos{\left(\frac{(2l-1)m \pi}{n}\right)} \ln \left\{\frac{\Gamma\left(1-\frac{2l-2a-1}{2n}\right)}{\Gamma(\frac{2l+2a-1}{2n})}\right\}, \ m+n \ \text{is even}
       \end{cases} \\
\end{align*}
where $b = \frac{\pi m}{n}$.

\paragraph{Medina \& Moll}

Using Theorem \ref{th:vardi} Medina \textit{et. al} \cite{moll:2008} calculated some known results found in \cite{ryzhik, vardi:1988, Adamchik:1997}. For example, by taking different $Q(x)$, they obtained the following integrals:
\allowdisplaybreaks{\begin{align*}
    & Q(x) = - \frac{\ln(1-x)}{x} = \sum\limits_{k=0}^{\infty} \frac{x^k}{k+1}, \quad  I_Q = \int_{0}^{1} \frac{\ln(1-x)}{x} \ln{\ln{\left(\frac{1}{x}\right)}} \ dx = \frac{\gamma \pi^2}{6} - \zeta(2)   \\
    & Q(x) =  \frac{\ln(1+x)}{x} = \sum\limits_{k=0}^{\infty} \frac{(-1)^kx^k}{k+1}, \quad  I_Q = \int_{0}^{1} \frac{\ln(1+x)}{x} \ln{\ln{\left(\frac{1}{x}\right)}} \ dx = \frac{ \pi^2}{12} (\ln{2}-\gamma) + \frac{1}{2} \zeta(2)   \\
    & Q(x) =  \frac{\text{PolyLog} [c,x]}{x} = \sum\limits_{k=0}^{\infty} \frac{x^k}{(k+1)^c}, \quad  I_Q = \int_{0}^{1} \frac{\text{PolyLog} [c,x]}{x} \ln{\ln{\left(\frac{1}{x}\right)}} \ dx = -\gamma \zeta(c+1) + \zeta^{'}(c+1)     \\
    & Q(x) =  -\frac{\text{Li}_{c-1}(-x)}{x} = \sum\limits_{k=0}^{\infty} \frac{(-1)^kx^k}{(k+1)^{c-1}}, \quad  I_Q = \int_{0}^{1} \frac{\text{Li}_{c-1}(-x)}{x} \ln{\ln{\left(\frac{1}{x}\right)}} \ dx \\
    & \qquad \qquad \qquad \qquad \qquad \qquad \qquad \qquad \qquad \qquad \qquad= \alpha [\gamma \zeta(c)-\zeta^{'}(c)] +(1+\alpha) \ln{2}\zeta(c) 
\end{align*}}
where $\alpha = 1-2^{1-c}$. Using Theorem \ref{th:vardi}, they also evaluated the following complicated integrals:
\begin{align*}
  R_{m,j}(a) = \int_{0}^{1} \frac{x^j  \ln{\ln{\left(\frac{1}{x}\right)}} \ dx}{(x+a)^{m+1}} , \  C_{m,j}(a,b) & = \int_{0}^{1} \frac{x^j  \ln{\ln{\left(\frac{1}{x}\right)}}  \ dx}{(x^2+ax+b)^{m+1}} \\
  & = D_{m,j}(r,\theta)  = \int_{0}^{1} \frac{x^j  \ln{\ln{\left(\frac{1}{x}\right)}} \ dx}{(x^2-2rx \cos{\theta}+r^2)^{m+1}}
\end{align*}

With $T_m(x)$ defined in \eqref{eq:T1}, they calculated the following integral for $1 \leq m \in \N$:
\begin{align*}
    \int_{0}^{1} \frac{T_{m-1}(x) \ln{\ln{\left(\frac{1}{x}\right)}}}{(x+1)^{m+1}} \ dx = (1-2^m) \ \zeta^{'}(1-m) + \left[\gamma(2^m-1)+2^m \ln{2}\right] \ \zeta(1-m)
\end{align*}

Using the above definitions they calculated the following generalized integral for any $m \in \N$ (Corollary 5.3,  \cite{moll:2008}):
\begin{align*}
    R_{m,0}(a) = \int_{0}^{1} \frac{  \ln{\ln{\left(\frac{1}{x}\right)}} }{(x+a)^{m+1}} \ dx = -\frac{\gamma}{m a^m(1+a)} - & \frac{\gamma}{a^{m+1}m!} \sum \limits_{j=2}^{m} \frac{S(m,j) T_{j-2}\left(\frac{1}{a}\right)}{\left(1+\frac{1}{a}\right)^j} \\
    & - \frac{1}{a^m m!} \sum \limits_{j=1}^{m} S(m,j) \ \text{Li}_{1-j} \left(-\frac{1}{a}\right)
\end{align*}
where, $S(m,j)$ are the unsigned Sterling number of first kind. The  result can be generalized via a recurrence relation for any $r \leq m$ given below (Theorem 6.1 in \cite{moll:2008}):
\begin{align*}
    R_{m,0}(a) = \sum\limits_{j=0}^{r} (-1)^j \binom{r}{j} a^{-r} R_{m-r+j,j}(a)
\end{align*}

For example with $m,j =0$ and $0 < \theta < 2\pi$, they derived the following integral
\begin{align*}
   D_{0,0}(1,\theta)  = \int_{0}^{1} \frac{ \ \ln{\ln{\left(\frac{1}{x}\right)}}}{(x^2-2rx \cos{\theta}+r^2)} \ dx = \frac{\pi}{2 \sin{\theta}} \left[\left(1-\frac{\theta}{\pi}\right)\ln{2\pi}+\ln{\left(\frac{\Gamma(1-\theta/2\pi)}{\Gamma(\theta/2\pi)}\right)}\right]
\end{align*}
Also the generalization for any $r \neq 1$, is derived as
\begin{align}
\label{eq:D1}
   D_{0,0}(r,\theta)  =  -\frac{\gamma}{r \sin{\theta}} \tan^{-1} \left(\frac{\sin{\theta}}{r-\cos{\theta}}
   \right) + \frac{1}{2 r i \sin{\theta}} \left[\text{Li}_{1}^{(1)}\left(\frac{e^{i \theta}}{r}\right)-\text{Li}_{1}^{(1)}\left(\frac{e^{-i \theta}}{r}\right)\right]
\end{align}
By choosing different values for the parameter $r, \theta$, (i.e., $r = 1, \theta = \frac{\pi}{2}, \frac{\pi}{3}, \frac{2\pi}{3}$) in equation \ref{eq:D1}, one can obtain most of the integrals provided in \cite{vardi:1988,ryzhik}.
And similarly, they derived the following generalization
\begin{align}
\label{eq:D2}
   D_{0,1}(r,\theta)  = -\frac{\gamma}{2} & \ln \left[\frac{r^2-2r \cos{\theta}+1}{r^2}\right] -\gamma \cot{\theta} \tan^{-1} \left(\frac{\sin{\theta}}{r-\cos{\theta}}
   \right) \nonumber \\
   & + \frac{1}{2 r i \sin{\theta}} \left[\Phi^{(1)}\left(\frac{e^{i \theta}}{r},1,1\right)-\Phi^{(1)}\left(\frac{e^{-i \theta}}{r},1,1\right)\right]
\end{align}
Now, using the above formulas one can derive closed form expressions for the generalized integrals $D_{m,j}(r,\theta)$
\begin{align*}
    D_{m,j}(r,\theta) = -\frac{1}{2rm \sin{\theta}} \frac{\partial}{\partial \theta} D_{m-1,j-1}(r,\theta) = \frac{r}{\cos{\theta}} D_{m,j-1}(r,\theta) +\frac{1}{2m \cos{\theta}} \frac{\partial}{\partial r} D_{m-1,j-1}(r,\theta)
\end{align*}

And finally, by using the same method authors recovered some of the Propositions provided in \cite{Adamchik:1997}

\paragraph{GR Tables revisited}
In their work \cite{moll:12}, Moll \textit{et. al} considered the following integral evaluation
\begin{align*}
    \int_{0}^{1} \frac{\ln{x}}{x^2-2ax+1} \ dx = \frac{1}{r_2-r_1} \left[\text{Li}_{2}\left(\frac{-1}{r_1}\right)-\text{Li}_{2}\left(\frac{-1}{r_2}\right)\right]
\end{align*}
where, $|a| < 1$ and $r_1, r_2$ are the imaginary roots of the denominator. By using the same approach they derived the following well known integrals
\begin{align*}
 \int_{0}^{1} \frac{x\ln{x}}{x^2-x+1} \ dx  = \frac{5 \pi^2}{36} - \frac{1}{6} \psi^{(1)} \left(\frac{1}{3}\right), \quad \int_{0}^{1} \frac{x\ln{x}}{x^2+x+1} \ dx  = -\frac{7 \pi^2}{54} + \frac{1}{9} \psi^{(1)} \left(\frac{1}{3}\right)
\end{align*}
They also calculated the following integrals involving higher powers of logarithms:
\begin{align*}
    \int_{0}^{1} \frac{1-x}{1-x^6} \ln^{2}x \ dx = \frac{2^3 \sqrt{3} \pi^3 + 351 \zeta(3)}{486}, \quad \int_{0}^{1} \frac{1-x}{1-x^6} \ln^{6}x \ dx = \frac{2^7 \times 7\sqrt{3} \pi^7 + 7 \times 1327995\zeta(7)}{26244}
\end{align*}
By looking at these integrals, one can rightly conjecture that the generalized integrals of these form has the following evaluation
\begin{align*}
     \int_{0}^{1} \frac{1-x}{1-x^6} \ln^{m-1}x \ dx = \frac{2^{m} \times a_m\sqrt{3} \pi^m + b_m \zeta(m)}{c_m}
\end{align*}
For a more detailed analysis for the computation of the above mentioned higher order integrals, see the work of Coffey \cite{coffey:2008}. Defining the integrals in mellin transformations author proved the following identity for $n \geq 0$ and $ p\geq 3$ odd integer:
\begin{align*}
    \int_{0}^{1} \frac{x^{p-3} \ln^n x}{1-x+x^2-...+x^{p-1}} & \ dx \\
    = & \frac{1}{(2p)^{n+1}} \left[\psi^{(n)}\left(1-\frac{1}{2p}\right)+\psi^{(n)}\left(1-\frac{1}{p}\right)-\psi^{(n)}\left(\frac{1}{2}-\frac{1}{2p}\right)-\psi^{(n)}\left(\frac{1}{2}-\frac{1}{p}\right)\right]
\end{align*}
For $\Re(s) > 0$ the author proved the following interesting identities
\begin{align*}
    & \int_{0}^{1} \frac{x^{p-3} \ln^s x}{1-x+x^2-...-x^{p-1}} \ dx = \frac{(-1)^s \Gamma(s+1)}{(2p)^{n+1}} \left[\zeta \left(s+1,1-\frac{1}{p}\right)+\zeta \left(s+1,1-\frac{2}{p}\right)\right], \quad p\geq 4 \ \text{even integer} \\
    & \int_{0}^{1} \frac{x^{p-3} \ln^s x}{1+x+x^2+...+x^{p-1}} \ dx = \frac{(-1)^s \Gamma(s+1)}{(2p)^{n+1}} \left[\zeta \left(s+1,1-\frac{1}{p}\right)+\zeta \left(s+1,1-\frac{2}{p}\right)\right], \quad p\geq 3 \ \text{ integer}
\end{align*}

\paragraph{GR Tables revisited} In \cite{moll:2007:27}, Medina \textit{et. al} considered the following types of integrals:
\begin{align*}
    & \int_{0}^{\infty} \frac{x^b-x^c}{1-x^a} \ dx = - \frac{\pi^2}{a^2} \frac{\sin(z-y)\sin(y+z)}{\sin^2(y)\sin^2(z)}, \ y = \frac{\pi b}{a}, z = \frac{\pi c}{a}, \quad I_n = \int_{0}^{1} \frac{x^{2n} \ln{x}}{(1+x^2)(1+x^4)^n} \ dx, n \in \N \\
    & \int_{0}^{1} \frac{x^{m-1} \ln^{a-1}(x)}{1+x} \ dx = (-1)^{a+m} \Gamma(a) \left[\frac{2^a-1}{2^{a-1}}\zeta(a) - \sum\limits_{k=1}^{m-1} \frac{(-1)^k}{k^a}\right], \quad \int_{0}^{\infty} \frac{x^{p-1} \ln{x}}{1-x^{q}} \ dx = -\frac{\pi^2}{q^2 \sin^2\left(\frac{\pi p}{q}\right)}  \\
    & f_n(a)  = \int_{0}^{1} \frac{x^a \ln^n{x}}{\sqrt{1-x^2}} \ dx = \lim_{s \rightarrow{a}} \left(\frac{\partial}{\partial s}\right)^n h(s), \quad \quad   h(s) = \int_{0}^{\frac{\pi}{2}} \sin^{s}t \ dt = \frac{1}{\sqrt{2}}\frac{\Gamma\left(\frac{s+1}{2}\right)}{\Gamma\left(\frac{s+2}{2}\right)}, \\
    &  \int_{0}^{1} \frac{x^{m-1} \ln^{a-1}(x)}{1-x} \ dx = (-1)^{a-1} \Gamma(a) \left[\zeta(a) - \sum\limits_{k=1}^{m-1} \frac{1}{k^a}\right] \\
    & \int_{0}^{\infty} \frac{x^{m-1} \ln{x}}{a-x} \ dx = \pi a^{m-1} \left[\ln{a} \cot{(\pi m)}-\frac{\pi}{\sin^(\pi m)}\right] \\
    & f_2(a) = \frac{\sqrt{\pi}}{8} \frac{\Gamma\left(\frac{a+1}{2}\right)}{\Gamma\left(\frac{a+2}{2}\right)} \left[\left\{\psi\left(\frac{a}{2}+1\right)-\psi\left(\frac{a+1}{2}\right)\right\}^2+\psi^{'}\left(\frac{a+1}{2}\right)-\psi^{'}\left(\frac{a}{2}+1\right)\right]
\end{align*}

They also calculated the following integral:
\begin{align*}
    \int_{0}^{1} \frac{x^{2n+1} }{\sqrt{1-x^2}} \ln^{2}(x) \ dx = -\frac{(2n)!!}{(2n+1)!!} \left[\frac{\pi^2}{12}+\sum\limits_{k=1}^{2n+1} \frac{(-1)^k}{k^2}- \left(\ln{2} + \sum\limits_{k=1}^{2n+1} \frac{(-1)^k}{k} \right)^2\right]
\end{align*}

In \cite{moll:10}, the following integrals are established using digamma function:
\begin{align*}
    \int_{0}^{1} \frac{\ln{x}}{\sqrt{x(1-x^2)}} \ dx = -\frac{\sqrt{2 \pi}}{8} \Gamma^2\left(\frac{1}{4}\right), \quad \int_{0}^{1} \ln{x} \sqrt{(1-x^2)^{2n-1}} \ dx = -\frac{(2n-1)!!}{4(2n)!!} \pi \left[\gamma + \ln{4}+\psi(n+1)\right]
\end{align*}

\paragraph{Blagouchine's Method} \footnote{this section is under construction}

There is a detailed review about Malmsten type integrals in his work \cite{blagouchine:2014}. For $a \geq 0$ and $\Re (b) >0$
\begin{align*}
    \int_{0}^{\infty} & \frac{\ln{(a^2+x^2)}}{\cosh^n{(bx)}}  \ dx \\
    = & \begin{cases}
       \frac{\pi A_n}{b} \ln{\left(  \frac{\Gamma \left(\frac{3}{4}+\frac{ab}{2\pi}\right)}{\Gamma \left(\frac{1}{4}+\frac{ab}{2\pi}\right)} \sqrt{\frac{2\pi}{b}}\right)} + \frac{1}{b} \sum\limits_{l=1}^{\frac{1}{2}(n-1)} \frac{D_{n,l}}{\pi^{2l-1}} \left\{\psi^{2l-1}\left(\frac{3}{4}+\frac{ab}{2\pi}\right)-\psi^{2l-1}\left(\frac{1}{4}+\frac{ab}{2\pi}\right)\right\}    , \ \text{for odd} \ n \\
       \frac{A_n}{b} \left\{\ln{\frac{\pi}{b}}+\psi\left(\frac{1}{2}+\frac{ab}{\pi}\right)\right\} + \frac{1}{b} \sum\limits_{l=1}^{\frac{1}{2}n-1} \frac{D_{n,l}}{\pi^{2l-1}} \ \psi^{2l-1}\left(\frac{1}{2}+\frac{ab}{\pi}\right), \  \text{for even} \ n
    \end{cases}
\end{align*}
where $A_n$ are the coefficients of the Maclaurin expansion of $2 (1-x)^{-\frac{1}{2}}$ for odd $n$ and for even $n$, $A_n^{-1}$ are equal to the coefficients in the Maclaurin expansion of $\frac{1}{2} (1-x)^{-\frac{3}{2}}$. There are no closed form formula available for the rational coefficients $D_{n,l}$ (see Table 1 in \cite{blagouchine:2014}). The following integrals were calculated in their work:
For $p = \frac{\pi m}{n}$, where $\Re (b) > 0$, and $m<n$ are positive integers and for any $a \geq 0$
\begin{multline*}
  \int_{0}^{\infty}  \frac{\cosh{(px)}\ \ln{(a^2+x^2)}}{\cosh{(bx)} + \cos{\Phi}} \ dx = \frac{2\pi}{b} \sin{\frac{m\Phi}{n}} \csc \Phi \csc \frac{m\pi}{n} \ln{\frac{2 \pi n}{b}} \\
  + \frac{2 \pi}{b \sin{\Phi}} \sum \limits_{l=0}^{n-1} [\cos{\frac{(2l+1)m\pi+m\Phi}{n}} \ln{\Gamma\left(\frac{2l \pi+\pi+ab+\Phi}{2\pi n}\right)} \\
  - \cos{\frac{(2l+1)m\pi-m\Phi}{n}} \ln{\Gamma\left(\frac{2l \pi+\pi+ab-\Phi}{2\pi n}\right)}]
\end{multline*}

The following generalized integrals were calculated: 
\begin{align*}
    \int_{0}^{\infty}  \frac{\cosh{(px)}\ \ln{(a^2+x^2)}}{\cosh^r{(bx)} + \sin^v{\Phi}} \ dx,  \qquad \int_{0}^{\infty}  \frac{\sinh{(px)}\ \ln{(a^2+x^2)}}{\sinh^r{(bx)} + \cos^v{\Phi}} \ dx 
\end{align*}
for combination of values $r = 1,2,3,4,5,6$ and $v = 0,1,2$. Furthermore, author also considered the following types of integrals:
\begin{align*}
    & \int_{0}^{1} \frac{x^{\frac{\alpha n}{2}-1}\ \ln{\ln{\frac{1}{x}}} \ dx}{1+x^{\alpha}+x^{2\alpha}+...+x^{n \alpha}} = \int_{1}^{\infty} \frac{x^{\frac{\alpha n}{2}-1}\ \ln{\ln{x}} \ dx}{1+x^{\alpha}+x^{2\alpha}+...+x^{n \alpha}} = \frac{2}{\alpha} \left\{Y_{n+1}-\frac{\pi}{2n+2} \tan \frac{\pi}{2n+2} \ \ln{\frac{\alpha}{2}}\right\} \\
    & \int_{0}^{1} \frac{x^{\frac{\alpha m}{2}-1}\ \ln{\ln{\frac{1}{x}}} \ dx}{1-x^{\alpha}+x^{2\alpha}-...+x^{m \alpha}} = \int_{1}^{\infty} \frac{x^{\frac{\alpha m}{2}-1}\ \ln{\ln{x}} \ dx}{1-x^{\alpha}+x^{2\alpha}-...+x^{m \alpha}} = \frac{2}{\alpha} \left\{X_{m+1}-\frac{\pi}{2m+2} \sec \frac{\pi}{2m+2} \ \ln{\frac{\alpha}{2}}\right\} 
\end{align*}
where $n = 1,2,3,...$ and $m = 2,4,6,...$, the function $Y_n$ and $X_n$ are defined in equation \eqref{eq:malm1} and \eqref{eq:malm2} respectively.

\paragraph{Moll's work} In their recent work, Moll et. al considered a wide range of integrals from the list of Gradshteyn \& Ryzhik. They developed generalized formulas for almost all of the results provided in \cite{ryzhik:7}. In the following paragraph, we will discuss some of these generalized formulas. In \cite{moll:2007:1}, authors considered the following type of integrals:
\begin{align}
    \int_{0}^{\infty} & \frac{\ln^{n-1}{x} \ dx}{(x-1)(x+a)}  = \frac{(-1)^n (n-1)!}{1+a} \left[\left\{1-(-1)^{n-1}\right\}\zeta(n)-\text{Li}_n\left(\frac{-1}{a}\right)+(-1)^{n-1} \text{Li}_n(-a)\right] \nonumber \\
    & = \frac{(-1)^n (n-1)!}{1+a} \left[\left\{1-(-1)^{n-1}\right\}\zeta(n)+ \frac{1}{n(1+a)} \sum\limits_{j=0}^{\floor*{\frac{n}{2}}} \binom{n}{2j} \left(2^{1-2j}-1\right)\pi_{2j} B_{2j} \left(\ln{a}\right)^{n-2j} \right]
\end{align}

In \cite{moll:2007:2}, authors considered the following types of integrals:
\begin{align}
  \int_{0}^{b} \frac{\ln{x}}{(x+r)^n} \ dx = \frac{\ln{r}}{(n-1) r^{n-1}} \left[\frac{(b+r)^{n-1}-r^{r-1}}{(b+r)^{n-1}}\right] + \frac{1}{r^{n-1}} h_{n} \left(\frac{b}{r}\right), \quad  h_n(b) = \int_{0}^{b} \frac{\ln{t}}{(1+t)^n} \ dt
\end{align}
The closed from expression for function $h_n(b)$ can be written as follows:
\begin{align*}
    h_n{b} = \frac{(1+b)^{n-1}-1}{(n-1)(1+b)^{n-1}} \ln{b} - \frac{(1+b)^{n-1}}{(n-1)(1+b)^{n-1}} \ln{(1+b)} + \frac{1}{(1+b)^{n-1}} Z_n(b)
\end{align*}
where the polynomials $Z_n(b)$ satisfy the following recurrence relation 
\begin{align*}
  Z_n(b) = \frac{(1+b)(n-2)}{n-1} Z_{n-1}(b) - \frac{(1+b)\left[(1+b)^{n-2}-1\right]}{(n-2)(n-1)}  
\end{align*}

\paragraph{Sofo} In his work \cite{sofo:2017}, Sofo considered the following integrals containing Gauss Hyper geometric and logarithmic functions. 
\begin{align*}
   I(m,p,q,t) = \frac{(-1)^{m-1}}{(m-1)!}\int_{0}^{1}  \pFq{2}{1}{1,\frac{1}{q}}{1+\frac{1}{q}}{x^q} \frac{1-x^p}{1-x} x^{-\frac{1}{t}}\ln^{m-1} x \ dx
\end{align*}

The following Theorem proved in \cite{sofo:2017}
\begin{theorem}
\label{th:sofo1}
(Theorem 1, \cite{sofo:2017}) Let $m, p \in \N, q \in \R \setminus \{-1,0\}$ and $t \in \R \setminus \{0\}$, then for $qkt-q-pt \neq 0$
\begin{align*}
    I(m,p,q,t) & = \frac{1}{q} \sum \limits_{k=1}^{p} \left(\frac{m}{qkt-q-pt}\right)^m \left(H_{\frac{kt-1}{pt}-1}-H_{\frac{1}{q}-1}\right) \\
    & + \frac{1}{qp^m} \sum \limits_{k=1}^{p} \sum \limits_{j=2}^{m} \left(\frac{pqt}{qkt-q-pt}\right)^{m+1-j} \left(H^{(j)}_{\frac{kt-1}{pt}-1}-\zeta(j)\right)
\end{align*}
\end{theorem}
Using special values for the parameters $m,p,q,t$ one can derive the following interesting integral:
\begin{align*}
    I(3,4,2,2) & = \frac{1}{2} \int_{0}^{1} \frac{(1-x^4)\ln^2 x}{(1-x)x^{\frac{5}{2}}}  \ln{\left(\frac{1+x^2}{1-x^2}\right)} \ dx \\
    & = \frac{16}{27} \left(14 \sqrt{2}-13\right) \pi + \frac{16}{3}\left(2 \sqrt{2}-5\right) \zeta(2) + \frac{\pi^3}{3} \left(3 \sqrt{2}-2\right) 
\end{align*}

He also proved the following Theorem:
\begin{theorem}
\label{th:sofo2}
(Theorem 2, \cite{sofo:2017}) Let $m, p, \mu \in \N, q \in \R \setminus \{-1,0\}$ and $t \in \R \setminus \{0\}$, then for $qkt-q-pt \neq 0$
\begin{align*}
   J(m,p,q,t,\mu) & = \frac{(-1)^{m-1}}{(m-1)!}\int_{0}^{1}   \frac{x^{p-\frac{1}{t}}(1-x^p)}{(1-x)q^{\mu+1}} \ln^{m-1} x \times \sum\limits_{r=0}^{\mu} \frac{(-1)^r}{q^r} \binom{\mu}{r} \Phi \left(x^p,1+r,1+\frac{1}{q}\right)    \ dx \\
   & = \sum \limits_{k=1}^{\infty} \frac{t^m\left[\frac{k-1}{p}\right]^\mu}{(kt-1)^m\left(q\left[\frac{k-1}{p}\right]+1\right)^{\mu+1}} = \sum \limits_{j=1}^{p} \sum\limits_{k=1}^{\infty} \frac{k^\mu}{\left(pk+j-\frac{1}{t}\right)^m(qk+1)^{\mu+1}}
\end{align*}
\end{theorem}

In \cite{sofo:2015}, the author derived the following integral identities (Lemma 6, \cite{sofo:2015}):
\begin{align*}
    \frac{1}{(2m-1)!} \int_{0}^{1} \frac{x \ \ln^{2m-1} x}{1-x} \Phi(-x,1,1+r) \ dx = M(m,r) -\zeta(2m) \Phi(-1,1,1+r)
\end{align*}
where $r \in \N \setminus \{0,1\}$ and $M(m,r)$ defined as follows:
\begin{align*}
    M(m,r) = \sum\limits_{k=1}^{\infty} \frac{(-1)^{k+1}}{k+r} H_{k}^{(2m)}
\end{align*}
In their work, author derived the following closed form expression for $M(m,r)$ (Lemma 4, \cite{sofo:2015}):
\begin{align*}
    M(m,r) = (-1)^{r+1} M(m,1) - (-1)^r \sum\limits_{k=1}^{r-1} \frac{H_k-H_{\left[\frac{k}{2}\right]}-\ln{2} -(-1)^{k+1}\ln{2}}{k^{2m}} + (-1)^r \sum\limits_{k=1}^{r-1} \sum\limits_{j=2}^{2m} \frac{(-1)^{k+1+j}}{k^{2m+1-j}} \bar{\zeta}(j)
\end{align*}
where $M(m,1)$ and $\bar{\zeta}(j)$ defined as follow:
\begin{align*}
    & \bar{\zeta}(z) = \sum\limits_{k=1}^{\infty} \frac{(-1)^{k+1}}{k^z} = \left(1-2^{1-z}\right) \zeta(z), \quad \bar{\zeta}(1) = \ln{2} \\
    & M(m,1) = \sum\limits_{k=1}^{\infty} \frac{(-1)^{k+1}H_k^{(2k)}}{k+1} = \bar{\zeta}(2m+1) - m \zeta(2m+1) + \ln{2} \ \bar{\zeta}(2m) + \sum\limits_{j=1}^{m-1} \bar{\zeta}(2j+1) \bar{\zeta}(2m-2j)
\end{align*}
There are several other complicated integrals can be found as follows:
\begin{theorem}
(Lemma 7, \cite{sofo:2015})
\label{th:sofo3}
\begin{align*}
    \frac{1}{(2m-1)!} \int_{0}^{1} \frac{x^2 \ \ln^{2m-1} x}{2(1-x)} \left[\Phi(x^2,1,\frac{1+r}{2})+ \Phi(x^2,1,\frac{2+r}{2}) \right] \ dx = \Upsilon(m,r) + \frac{1}{2} \zeta(2m) \left(H_{\frac{r-1}{2}}-H_{\frac{r}{2}}\right)
\end{align*}
\end{theorem}

where $r \in \N \setminus \{0,1\}$ and for $r \geq 2, m \in \N$ the function $\Upsilon(m,r)$ is defined as follows:
\begin{align*}
    \Upsilon(m,r) = \sum\limits_{k=1}^{\infty} \frac{H_{2k}^{(2m)}}{(2k+r-1)(2k+r)}  =  M(m,r) - \frac{H_{\frac{r-1}{2}}}{2(r-1)^{2m}} + \sum\limits_{j=2}^{2m} \frac{(-1)^j \zeta(j)}{2^j (r-1)^{2m+1-j}} 
\end{align*}
They also derived the following integrals:
\begin{theorem}
(Theorem 3, \cite{sofo:2015}) For $m,k \in \N$, the following relations hold:
\label{th:sofo4}
\begin{align*}
    & \frac{1}{(1+k)(2m-1)!}\int_{0}^{1}  \pFq{2}{1}{1,1}{k+2}{-x} \frac{x \ \ln^{2m-1} x }{1-x}  \ dx  = X(m,k,1) - \frac{\zeta(2m)}{k+1} \ \pFq{2}{1}{1,1}{k+2}{-1}
\end{align*}
and 
\begin{align*}
    & \frac{1}{(1+k)(2m-1)!}\int_{0}^{1}  \pFq{2}{1}{1,2}{k+2}{-x} \frac{x \ \ln^{2m-1} x }{1-x}  \ dx  = X(m,k,0) - \frac{\zeta(2m)}{k+1} \ \pFq{2}{1}{1,2}{k+2}{-1} 
\end{align*}
\end{theorem}
where the function $X(m,k,p)$ for $p = 0,1$ is given by the following:
\begin{align*}
     & X(m,k,0) = \sum\limits_{n=1}^{\infty} \frac{(-1)^{n+1}H_{n}^{(2m)}}{\binom{n+k}{k}}  =   \sum\limits_{j=1}^{k} (-1)^{1+j} j \binom{k}{j} M(m,j) \\
     & X(m,k,1) = \sum\limits_{n=1}^{\infty} \frac{(-1)^{n+1}H_{n}^{(2m)}}{n\binom{n+k}{k}}  =  M(m,0)- \sum\limits_{j=1}^{k} (-1)^{1+j} j \binom{k}{j} M(m,j)
\end{align*}
By using  some elementary operations, the following interesting identity is derived:
\begin{align*}
    \frac{1}{(2m-1)!} \int_{0}^{1} \frac{\ln^{2m-1} x \ \text{Li}_p(-x)}{x} \ dx = \bar{\zeta}(2m+p)
\end{align*}
for $p \in \N$, which can be seen as a analogous identity of the following one:
\begin{align*}
     \int_{0}^{1} \frac{\ln^{m} x \ \text{Li}_p(x)}{x} \ dx = (-1)^m m!\zeta(m+p+1)
\end{align*}

In \cite{sofo:2018}, Sofo derived the following integral:

\begin{theorem}
(Lemma 1.1, \cite{sofo:2018}) Let $k$ be a positive integer, then
\label{th:sofo6}
\begin{align*}
    & X(k,0) = \sum\limits_{n=1}^{\infty} \frac{(-1)^{n+1} H_{kn}}{n} = - k \int_{0}^{1} \frac{x^{k-1}\ \ln{(1-x)}}{1+x^k} \ dx = \frac{1+k^2}{4k} \zeta(2) - \frac{1}{2} \sum\limits_{j=0}^{k-1} \ln^2 \left[2 \sin \left(\frac{2j+1}{2k} \pi\right)\right]
\end{align*}
\end{theorem}

\begin{theorem}
(Lemma 1.2, \cite{sofo:2018}) Let $k$ be a positive integer, then
\label{th:sofo7}
\begin{align*}
    & X(k,1) =\sum\limits_{n=1}^{\infty} \frac{(-1)^{n+1} H_{kn}}{n+1} =  - k \int_{0}^{1} \frac{\ln{(1-x)}}{x} \left[\frac{1}{1+x^k}-\frac{\ln{(1+x^k)}}{x^k}\right] \ dx \\
    & \qquad \qquad \qquad \qquad \qquad \qquad = -X(k,0)+ \frac{1}{2k} \zeta(2)- H_{k-1} \ln{2} - \frac{1}{2} \sum\limits_{j=1}^{k-1} \frac{1}{j} \left(H_{-\frac{1}{2k}}-H_{-\frac{k+j}{2k}}\right)
\end{align*}
\end{theorem}

And the generalization $X(k,r)$, where $k$ and $r$ are positive integers, is given by the following theorem:
\begin{theorem}
(Lemma 1.3, \cite{sofo:2018}) Let $k$ and $r$ be positive integers, then
\label{th:sofo7_1}
\begin{align*}
    & X(k,r) =   \sum\limits_{n=1}^{\infty} \frac{(-1)^{n+1} H_{kn}}{n+r} = - \frac{k}{1+r} \int_{0}^{1} x^{k-1} \ln{(1-x)} \pFq{2}{1}{2,1+r}{2+r}{-x^k} \ dx \\
    & = (-1)^{k+1} X(k,1) - \frac{(-1)^{r} \ln{2}}{k} \left(H_{\left[\frac{r-1}{2}\right]} - H_{r-1}\right) -\frac{(-1)^{r}}{2k} \sum\limits_{j=1}^{r-1} \frac{(-1)^j}{j} \left(H_{\frac{j-1}{2}}-H_{\frac{j}{2}}\right) \\
    & - \frac{(-1)^{r}}{2} \sum\limits_{j=1}^{r-1} (-1)^j \left(H_{\frac{j-1}{2}}-H_{\frac{j}{2}}\right) \left(H_{kj+k-1}-H_{kj}\right) - \frac{(-1)^{r}}{2} \sum\limits_{j=1}^{r-1} \sum\limits_{i}^{k-1} \frac{(-1)^j}{kj+i} \left(H_{-\frac{1}{2k}}-H_{-\frac{k+i}{2k}}\right)
\end{align*}
\end{theorem}

He also derived the following interesting identities.
\begin{theorem}
(Lemma 1.4, \cite{sofo:2018})
\label{th:sofo8}
Let $k$ and $r \geq 2$ be positive integers, then
\begin{align*}
    & Y(k,r) = \frac{k}{2} \int_{0}^{1} x^{2k-1} \ln{(1-x)} \left[(r-1) \Phi\left(x^{2k},1,\frac{r+1}{2}\right)- \Phi\left(x^{2k},1,\frac{r+2}{2}\right) \right] \ dx \\
    & = X(k,r) + \frac{H_{\frac{r-1}{2}}}{2} \left(H_{kr-1}-H_{kr-k-1}\right) - \frac{1}{2} \sum\limits_{j=1}^{k-1} \frac{1}{j+rk-k} H_{-\frac{1}{2k}} = \sum\limits_{n=1}^{\infty} \frac{H_{2kn}}{(2n+r-1)(2n+r)}
\end{align*}
\end{theorem}
Two special cases of the above Theorem ($r = 0, 1$) can be derived as follows:
\begin{align*}
    & \frac{k}{2} \int_{0}^{1} x^{k-1} \ln{(1-x)} \ln{\left(\frac{1+x^k}{1-x^k}\right)} \ dx = - X(k,0) -H_k \ln{2} + \frac{1}{2} \sum\limits_{j=1}^{k-1} \frac{1}{j-k} H_{-\frac{j}{2k}} \\
    & k \int_{0}^{1}  \frac{\ln{(1-x)}}{x}  \left[1- \frac{1}{2x^k} \ln{\left(\frac{1+x^k}{1-x^k}\right)} \right] \ dx = X(k,1) + \frac{1}{4k} \zeta(2) - \frac{1}{2} \sum\limits_{j=1}^{k-1} \frac{1}{j} H_{-\frac{j}{2k}}
\end{align*}

There are some  other interesting identities derived in \cite{sofo:2018}, for $k\geq 1$ and p are real positive integers,

\begin{theorem}
(Theorem 2.1, \cite{sofo:2018}) Let $k\geq 1$ and $p$ be real positive integers, then
\label{th:sofo9}
\begin{align*}
    & \Delta (k,p) = \sum\limits_{n=1}^{\infty} \frac{(-1)^n H_{kn}}{\binom{n+p}{p}} = \frac{k}{p+1} \int_{0}^{1} x^{k-1} \ln{(1-x)} \pFq{2}{1}{2,2}{2+p}{-x^k} \ dx = \sum\limits_{r=0}^{p} (-1)^r r \binom{p}{r} X(k,r)
\end{align*}
\end{theorem}

\begin{theorem}
(Theorem 2.2, \cite{sofo:2018}) Let $k\geq 1$ and $p$ be real positive integers, then
\label{th:sofo10}
\begin{align*}
    & M (k,p) = \sum\limits_{n=1}^{\infty} \frac{(-1)^n H_{kn}}{n\binom{n+p}{p}} = \frac{k}{p+1} \int_{0}^{1} x^{k-1} \ln{(1-x)} \pFq{2}{1}{1,2}{2+p}{-x^k} \ dx = \sum\limits_{r=0}^{p} (-1)^r  \binom{p}{r} X(k,r)
\end{align*}
\end{theorem}

\begin{theorem}
\label{th:sofo11}
(Theorem 2.3, \cite{sofo:2018})
For $p \in \N \cup \{0\}$ and $k \in \N$
\begin{align*}
    & \frac{1}{k(p+1)} \sum \limits_{r=1}^{k-1} \int_{0}^{1} x^{-\frac{r}{k}} \ln{(1-x)} \left(k \ \pFq{2}{1}{1,2}{2+p}{-x}- r \ \pFq{2}{1}{1,1}{2+p}{-x} \right) \ dx \\
    & = - k M(k,p) - \frac{k \ln{k}}{p+1} \pFq{2}{1}{1,1}{2+p}{-1} - S(p) = \sum \limits_{r=1}^{k-1} \sum\limits_{n=1}^{\infty} \frac{(-1)^{n+1}H_{n-\frac{r}{k}}}{n \binom{n+p}{p}}
\end{align*}
\end{theorem}
where $S(p)$ is given by
\begin{align*}
    S(p) = \frac{1}{2} \zeta(2)  + \left(2^{p-1}-1\right) \ln^2 2 + & \sum \limits_{m=1}^{p} m \binom{p}{m} \left[ \left(2H_{m-1}-H_{\left[\frac{m-1}{2}\right]}\right) \ln{2}+ H_{m-1} \left(H_{\left[\frac{m}{2}\right]}-H_m\right) \right] \\
    & - \sum \limits_{m=1}^{p} \sum\limits_{j=1}^{m-1}  \frac{(-1)^j \ m}{j} \binom{p}{m} \left(H_{\left[\frac{m-j}{2}\right]}-H_{m-j} + \frac{j H_j}{j+1}
    \right)
\end{align*}
where $[x]$ denotes the integer part of $x$.

\begin{theorem}
\label{th:sofo12}
(Theorem 2.4, \cite{sofo:2018}) Let $k$ and $p$ be real real positive integers, then
\begin{align*}
    & \Omega(k,p) = \frac{-4k}{(p+1)(p+2)} \int_{0}^{1} x^{2k-1} \ln{(1-x)}  \ \pFq{3}{2}{1,\frac{3}{2},2}{\frac{p+3}{2},\frac{p+4}{2}}{x^{2k}} \ dx =   \sum\limits_{n=1}^{\infty} \frac{H_{2kn}}{n \binom{2n+p}{p}} \\
    & = \frac{2}{p} \sum\limits_{r=1}^{p} (-1)^{1+r} r \binom{p}{r} X(k,r) + \frac{2}{p(p+1)} \sum\limits_{j=0}^{k-1} \frac{1}{2k-j} \ \pFq{4}{3}{\frac{1}{2},1,1,\frac{2k-j}{2k}}{\frac{p+2}{2},\frac{p+3}{2},\frac{4k-j}{2k}}{1}
\end{align*}
\end{theorem}

\begin{theorem}
\label{th:sofo13}
(Theorem 2.5, \cite{sofo:2018})
Let $k$ and $p$ be real real positive integers, then
\begin{align*}
    & \Xi(k,p) = \frac{-4k}{(p+1)(p+2)} \int_{0}^{1} x^{2k-1} \ln{(1-x)}  \ \pFq{3}{2}{\frac{1}{2},1,2}{\frac{p+3}{2},\frac{p+4}{2}}{x^{2k}} \ dx =   \sum\limits_{n=1}^{\infty} \frac{H_{2kn}}{n(2n-1) \binom{2n+p}{p}} \\
    & = \frac{2}{p+1} \sum\limits_{r=0}^{p} (-1)^{r} \binom{p}{r} X(k,r) + \frac{2}{(p+1)^2} \sum\limits_{j=0}^{k-1} \frac{1}{2k-j} \ \pFq{4}{3}{1,1,\frac{3}{2},\frac{2k-j}{2k}}{\frac{p+2}{2},\frac{p+3}{2},\frac{4k-j}{2k}}{1}
\end{align*}
\end{theorem}

\begin{theorem}
\label{th:sofo14}
(Theorem 2.6, \cite{sofo:2018})
Let $k$ and $p$ be real real positive integers, then
\begin{align*}
    & \Upsilon(k,p) = \frac{4k}{(p+1)(p+2)} \int_{0}^{1} x^{2k-1} \ln{(1-x)}  \ \pFq{3}{2}{\frac{1}{2},2,2}{\frac{p+3}{2},\frac{p+4}{2}}{x^{2k}} \ dx =   \sum\limits_{n=1}^{\infty} \frac{H_{2kn}}{(2n-1) \binom{2n+p}{p}} \\
    & = \frac{2}{(p+1)^2}  \sum\limits_{j=0}^{k-1} \frac{1}{2k-j} \ \pFq{4}{3}{\frac{1}{2},1,1,\frac{2k-j}{2k}}{\frac{p+2}{2},\frac{p+3}{2},\frac{4k-j}{2k}}{1}+ \frac{2}{p(p+1)} \sum\limits_{j=0}^{k-1} \frac{1}{2k-j} \ \pFq{4}{3}{1,1,\frac{3}{2},\frac{2k-j}{2k}}{\frac{p+2}{2},\frac{p+3}{2},\frac{4k-j}{2k}}{1}
\end{align*}
\end{theorem}

Using special values for the parameters ($k,r,p$) in the above mentioned Theorems \ref{th:sofo6}-\ref{th:sofo14}, one can obtain the following interesting and complicated integrals (Corollary 2.7 and 2.8 in \cite{sofo:2018}):
\begin{align*}
    & \int_{0}^{1} \frac{1+x^{2k}}{x^{k+1}} \ln{(1-x)} \ln{\left(\frac{1+x^k}{1-x^k}\right)} \ dx = -\frac{3+4k^2}{2k^2} \zeta(2) -\frac{2}{k^2} \ln{2} - \frac{4}{k} H_{k-1} + \frac{1}{k} \sum\limits_{j=1}^{k-1} \frac{1}{k-j} \left(H_{\frac{j}{2k}}-H_{-\frac{j}{2k}}\right) \\
    & \int_{0}^{1} \frac{1+x}{\sqrt{x}} \ln^2{x} \ln{\left(\frac{1+x}{1-x}\right)} \ dx = -\frac{10}{3} \zeta(2) -\frac{14}{3} \zeta(3)+\frac{26}{27} \pi +\frac{32}{9} G -\frac{56}{27} \ln{2}+\frac{\pi^3}{12} +\frac{104}{27} \\
    & \int_{0}^{1} \frac{1-x^4}{(1-x)x^{5/2}} \ln^2{x} \ln{\left(\frac{1+x^2}{1-x^2}\right)} \ dx = \frac{16}{27} \left(14 \sqrt{2}-13\right) \pi + \frac{16}{3} \left(2\sqrt{2}-5\right) \zeta(2) + \frac{\pi^3}{3} \left(3\sqrt{2}-2\right)
\end{align*}

\section{Conclusion}
In this work, we provided a systematic review of logarithmic integrals arising from analytic number theory and complicated physical models. This work will be continuously updated and reviewed with the recent advances in this field.

\bibliographystyle{unsrt}  
\bibliography{references}  

\begin{thebibliography}{10}

\bibitem{ryzhik}
Alan Jeffrey, Daniel Zwillinger, I.S. Gradshteyn, and I.M. Ryzhik.
\newblock {\em Table of Integrals, Series, and Products (Seventh Edition)}.
\newblock Academic Press, Boston, seventh edition edition, 2007.

\bibitem{vardi:1988}
Ilan Vardi.
\newblock Integrals, an introduction to analytic number theory.
\newblock {\em The American Mathematical Monthly}, 95(4):308--315, 1988.

\bibitem{Adamchik:1997}
Victor Adamchik.
\newblock A class of logarithmic integrals.
\newblock In {\em Proceedings of the 1997 International Symposium on Symbolic
  and Algebraic Computation}, ISSAC '97, pages 1--8, New York, NY, USA, 1997.
  ACM.

\bibitem{boros:2006}
G.~Boros and V.~Moll.
\newblock {\em Irresistible integrals: symbolics, analysis and experiments in
  the evaluation of integrals}.
\newblock Cambridge Univ. Press., 2006.

\bibitem{moll:2008}
Luis Medina and Victor Moll.
\newblock A class of logarithmic integrals, 2008.

\bibitem{ripon:2015}
Sarowar~Morshed Ripon.
\newblock Generalization of a class of logarithmic integrals.
\newblock {\em Integral Transforms and Special Functions}, 26(4):229--245,
  2015.

\bibitem{ripon:2014}
Sarowar~Morshed Ripon.
\newblock Generalization of harmonic sums involving inverse binomial
  coefficients.
\newblock {\em Integral Transforms and Special Functions}, 25(10):821--835,
  2014.

\bibitem{ripon:2016}
Sarowar~Morshed Ripon.
\newblock A generalized inverse binomial summation theorem and some
  hypergeometric transformation formulas.
\newblock {\em International Journal of Combinatorics}, page~14, 2016.

\bibitem{ryzhik:7}
Alan Jeffrey, Daniel Zwillinger, I.S. Gradshteyn, and I.M. Ryzhik.
\newblock 3–4 - definite integrals of elementary functions.
\newblock In {\em Table of Integrals, Series, and Products (Seventh Edition)},
  pages 247 -- 617. Academic Press, Boston, seventh edition edition, 2007.

\bibitem{blagouchine:2014}
Iaroslav~V. Blagouchine.
\newblock Rediscovery of malmsten's integrals, their evaluation by contour
  integration methods and some related results.
\newblock {\em The Ramanujan Journal}, 35(1):21--110, Oct 2014.

\bibitem{potts:1978}
R.~J. Baxter, H.~N.~V. Temperley, and S.~E. Ashley.
\newblock Triangular potts model at its transition temperature, and related
  models.
\newblock {\em Proceedings of the Royal Society of London. Series A,
  Mathematical and Physical Sciences}, 358(1695):535--559, 1978.

\bibitem{malmsten:1}
Almgren T.A. Camitz G. Danelius D. Moder D.H. Selander E. Grenander J.M.A.
  Themptander S. Trozelli L.M. Föräldrar Ä. Ossbahr G.E. Föräldrar D.H.
  Ossbahr C.O. Lindhagen C.A. Moder D.H. Syskon Ä. Lemke O.V. Fries C.
  Laurenius L. Leijer E. Gyllenberg G. Morfader M.V. Linderoth~A. Malmsten,
  C.J.
\newblock Specimen analyticum, theoremata quædam nova de integralibus
  definitis, summatione serierum earumque in alias series transformatione
  exhibens (eng. trans.: “some new theorems about the definite integral,
  summation of the series and their transformation into other series”).
\newblock {\em Upsaliæ, excudebant Regiæ academiæ typographi. Uppsala,
  Sweden}, pages Dissertation, in 8 parts, April–June 1842.

\bibitem{malmsten:2}
C.J. Malmstén.
\newblock De integralibus quibusdam definitis seriebusque infinitis (eng.
  trans.: “on some definite integrals and series”).
\newblock {\em J. Reine Angew. Math.}, 38:1--39, (1849) [work dated May 1,
  1846].

\bibitem{potts:1997}
Robert~M. Ziff, Steven~R. Finch, and Victor~S. Adamchik.
\newblock Universality of finite-size corrections to the number of critical
  percolation clusters.
\newblock {\em Phys. Rev. Lett.}, 79:3447--3450, Nov 1997.

\bibitem{moll:12}
Victor~H. Moll and Ronald~A. Posey.
\newblock The integrals in gradshteyn and ryzhik. part 12: Some logarithmic
  integrals, 2010.

\bibitem{coffey:2008}
M.~Coffey.
\newblock Evaluation of certain mellin transformations in terms of the trigamma
  and polygamma functions.
\newblock {\em Contemporay Mathematics. Special Functions and Orthogonal
  Polynomials}, 471:85--104, 2008.

\bibitem{moll:2007:27}
Luis~A. Medina and V.~Moll.
\newblock The integrals in gradshteyn and ryzhik. part 27: More logarithmic
  examples.
\newblock {\em Scientia}, 2014.

\bibitem{moll:10}
Luis~A. Medina and Victor~H. Moll.
\newblock The integrals in gradshteyn and ryzhik. part 10: the digamma
  function, 2007.

\bibitem{moll:2007:1}
Victor~H. Moll.
\newblock The integrals in gradshteyn and rhyzik. part 1: A family of
  logarithmic integrals, 2007.

\bibitem{moll:2007:2}
V.~Moll.
\newblock The integrals in gradshteyn and ryzhik. part 2: Elementary
  logarithmic integrals.
\newblock {\em Scientia}, 14:7--15, 2007.

\bibitem{sofo:2017}
Anthony Sofo.
\newblock A master integral in four parameters.
\newblock {\em Journal of Mathematical Analysis and Applications}, 448(1):81 --
  92, 2017.

\bibitem{sofo:2015}
Anthony Sofo.
\newblock Integrals of logarithmic and hypergeometric functions.
\newblock {\em Communications in Mathematics}, 24, 09 2015.

\bibitem{sofo:2018}
Anthony Sofo.
\newblock Evaluation of integrals with hypergeometric and logarithmic
  functions.
\newblock {\em Open Mathematics}, 16(1):63–74, 2018.

\end{thebibliography}


\end{document}